\documentclass[12pt,a4paper]{article}
\usepackage{url}
\usepackage[utf8]{inputenc}
\usepackage[english]{babel}
\usepackage{amsmath}
\usepackage{amsthm}
\usepackage{amsfonts}
\usepackage{caption}
\usepackage{xcolor}
\usepackage{subcaption}
\usepackage{amssymb}
\usepackage{multicol}
\usepackage{tcolorbox}
\usepackage{enumerate}
\usepackage{hyperref}
\usepackage{graphicx}
\usepackage{parskip}
\usepackage[all]{xy}
\usepackage[left=2.54cm,right=2.54cm,top=2.54cm,bottom=2.54cm]{geometry}

\newtheorem{theorem}{Theorem}[section]

\newtheorem{lemma}{Lemma}[section]
\newtheorem{proposition}{Proposition}[section]
\newtheorem{corollary}{Corollary}[section]

\author{}
\title{\textbf{On transverse-universality of twist knots}\\}
\author{Sebastian Zapata}
\date{}

\begin{document}
\maketitle

\begin{abstract}
In the search for transverse-universal knots in the standard contact structure on $\mathbb{S}^3$, we present a classification of the transverse twist knots with maximal self-linking numbers, that admit only overtwisted contact branched coverings.  As a direct consequence, we obtain an infinite family of transverse knots in $(\mathbb{S}^3,\xi_{std})$ that are not transverse-universal, although they are universal in the topological sense.

\end{abstract}

\section{Introduction}

In 2017, Roger Casals and John Etnyre extended the classical existence theorem of universal links in $\mathbb{S}^3$ to the category of contact $3$-manifolds. They proved that every contact  $3$-manifold is the contact branched covering over $(\mathbb{S}^3,\xi_{std})$, with branch locus a fixed four-component transverse link \cite{CaEt}. Later, Jesús Rodríguez proved that given a transverse-universal link in $(\mathbb{S}^3,\xi_{std})$, it is possible to reduce the number of its components, maintaining the transverse-universality property. In particular, this argument guaranteed the existence of transverse-universal knots in the standard contact structure of $\mathbb{S}^3$ \cite{Ro}.

Nevertheless, the proof given by Jesús Rodr\'iguez used an algorithmic process to reduce the number of components by rapidly increasing the number of crossings in consecutive steps. It follows that the search of easily manipulable transverse-universal knots is still an open question, and a lot of the classical universal knots are good candidates for this. 

In her doctoral thesis, Meredith Casey proved in 2013 that the figure eight knot  is not transverse-universal, even if it is universal in the topological sense. She showed that every contact branched covering over $(\mathbb{S}^3,\xi_{std})$ with branch locus that knot must be overtwisted, and hence,  we cannot obtain any tight contact structure as a covering branched over it {\cite[Theorem~4.3.1]{Cas}}. We will show a  direct extension of her work by working with the complete family of twist knots, discarding a subfamily in the same way as Casey, and proving that the other topological universal members have associated tight contact branched coverings. More precisely, we prove the following theorem:

\begin{theorem}\label{Maintheo}
 If $m\geq2$ or $m=2k-1\leq-3$, then every contact branched covering with branch locus a transverse representative of the twist knot $K_m$ is overtwisted. 
\end{theorem}

The case studied by Casey corresponds to $m=2$, and by analyzing the other possibilities for $m$, we are going to see that for every $m\leq -4$ even, $K_m$ admits tight contact branched coverings, and therefore, $K_{-4}$ will be the simplest knot in this family that we cannot discard by using the above method.

\section{Background}

 For this part, we will review the concepts we need to prove the Theorem \ref{Maintheo}. Except for $\mathbb{R}^3$, all our $3$-manifolds will be closed and oriented.

\subsection{Branched coverings}

Recall that a \textit{branched covering} between $3$-manifolds is a map $\pi:M\rightarrow N$ for which exists a link $L\subseteq N$ called branch locus, such that the restriction $M-\pi^{-1}(L)\rightarrow N-L$ is a finite covering, and $ \pi^{-1}(L)$ is contained in a neighborhood in which this map can be viewed as the function $\mathbb{S}^1\times\mathbb{D}^2\rightarrow\mathbb{S}^1\times\mathbb{D}^2$: $(\theta, z)\rightarrow(n\theta, z^m)$ for certain positive integers $n,m$.

Consider now an ordinary $n$-sheeted covering map $\pi:M\rightarrow N$ with $M$ connected, and fix $x_0\in N$ and a labeling $\pi^{-1}(x_0)=\{x_1,\dots,x_n\}$. Then every loop $\gamma\subseteq N$ based at $x_0$, lifts to a unique $\tilde\gamma_i\subseteq M$ with starting point $x_i$, that ends in a certain point $x_j$. In this way, we can define a permutation $\sigma\in S_n$ with $\sigma(i)=j$. This process induces a map $m:\pi_1(N,x_0)\rightarrow S_n$ defined by $m([\gamma])=\sigma$ as before; it is called the \textit{monodromy map} associated with the covering. For branched coverings, the monodromy map is the monodromy of the associated ordinary covering.

Let $L$ be a link in a $3$-manifold $N$. We say that $L$ is \textit{universal} if, for every closed, oriented $3$-manifold $M$, there exists a finite branched covering $\pi:M\rightarrow N$ with branch locus the given link. See \cite{HiLoMo1} for more details.

\subsection{Contact 3-manifolds}

A \textit{co-orientable contact structure} on a $3$-manifold $M$ is a smooth distribution of hyperplanes $\xi$ defined as the kernel of a 1-form  $\alpha$ with $\alpha\wedge d\alpha\neq0$. As usual, we restrict our work to positive contact structures, that is, $\alpha\wedge d\alpha>0$. A link $L\subseteq M$ is \textit{Legendrian} if it is tangent to the distribution at every point and \textit{transverse} if it is transverse everywhere.

We say that a contact structure is \textit{overtwisted} if there exists an embedded disk $D$, whose boundary is tangent to the plane field distribution at every point. Otherwise, we call it \textit{tight}.

 By removing a point from $\mathbb{S}^3$, we will view $(\mathbb{R}^3,\xi_{std})\subseteq(\mathbb{S}^3,\xi_{std})$, where $\xi_{std}$ denotes the \textit{standard contact structure}, defined in $\mathbb{R}^3$ by $\xi_{std}=\ker(dz-ydx)$ in the usual Cartesian coordinates. This is equivalent to $(\mathbb{R}^3,\xi_{rot})$ with ${\xi_{rot}=\ker(dz+xdy-ydx)}$, but this last one will be useful for the treatment of transverse links and their relation with closed braids.
 
\begin{figure}
    \centering
    \includegraphics[width=11cm]{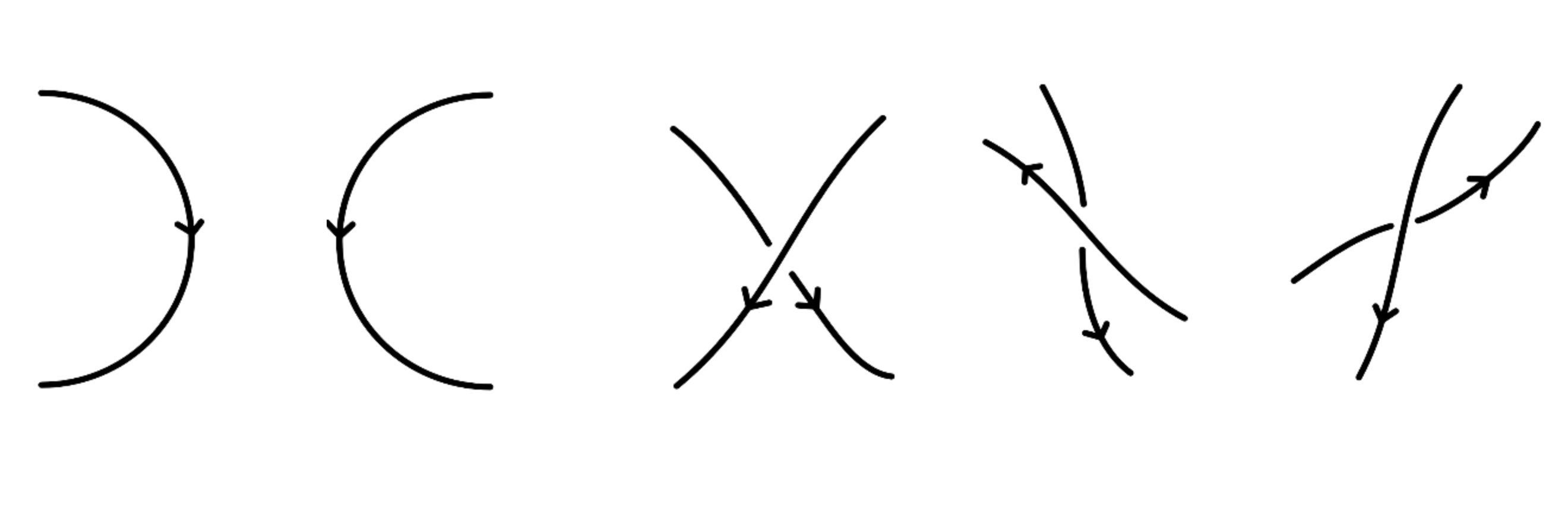}
    \caption{Forbidden portions and crosses for a front projection}
    \label{forbidden crosses}
\end{figure}

 If we take a knot $L\subseteq(\mathbb{R}^3,\xi_{std})$, its \textit{front projection} is defined by $\pi:\mathbb{R}^3\rightarrow\mathbb{R}^2$, $(x,y,z)\mapsto(x,z)$. The condition that $L$ is transverse to the contact planes implies the existence of forbidden portions and crosses in its front projection, represented in Figure \ref{forbidden crosses}. Conversely, every given projection of a link in the $xz$ plane that does not contain portions of the figure is the front projection of some transverse knot.

There is also an equivalent point of view for contact $3$-manifolds. An \textit{open book decomposition} of $M$ is an oriented link $L\subseteq M$ called \textit{the binding} and a fibration $\pi:(M-L)\rightarrow \mathbb{S}^1$ such that each fiber $\Sigma_\theta:=\pi^{-1}({\theta})$ is the interior of a compact surface $\Sigma$, called \textit{page}, with boundary $L$.  An alternative description, called \textit{abstract open book decomposition}, comes from an automorphism of the surface $\phi$, called \textit{monodromy}, and the construction of the mapping torus. These two notions are compatible in the sense that they recover, up to equivalence, diffeomorphic manifolds, see  \cite{Et} for further details. By the Giroux correspondence, there is a $1$ to $1$ equivalence between contact isotopy classes of oriented contact structures, and open book decompositions up to positive Hopf stabilization \cite{Et}.

The particular importance of this result is that for $(\mathbb{R}^3,\xi_{rot})$ we can choose an associated decomposition of our manifold, by taking the $z$-axis as the binding, half-planes with boundary the given axis as pages, and the identity as the surface automorphism. Through compactification, it induces an open book decomposition of $\mathbb{S}^3$ with page $D^2$ and monodromy the identity. With this, every transverse knot in this contact manifold can be braided in $\mathbb{R}^3$ around the $z$-axis \cite{Ben}. That is, each of those transverse knots is equivalent to a braid through isotopies that preserve the transverse property.

Let $\pi:M\rightarrow N$ be a finite branched covering with branch locus $L$. By fixing a contact structure $\xi_N=\ker(\alpha)$ on $N$ with $L$ transverse, we can define a unique contact structure $\xi_M$ on $M$ by considering an arbitrarily small deformation of $\pi^*\alpha$ near $\pi^{-1}(L)$ \cite{Go}. We denote this by $\pi:(M,\xi_M)\rightarrow (N,\xi_N)$ and call it \textit{contact branched covering}. In this way, we can say that a transverse link $L$ in $(N,\xi_N)$ is \textit{transverse-universal} if, for every closed, oriented $3$-manifold $M$ and every contact structure $\xi_M$ on it, there exists a contact branched covering as before.
\begin{figure}
    \centering
    \includegraphics[width=13cm]{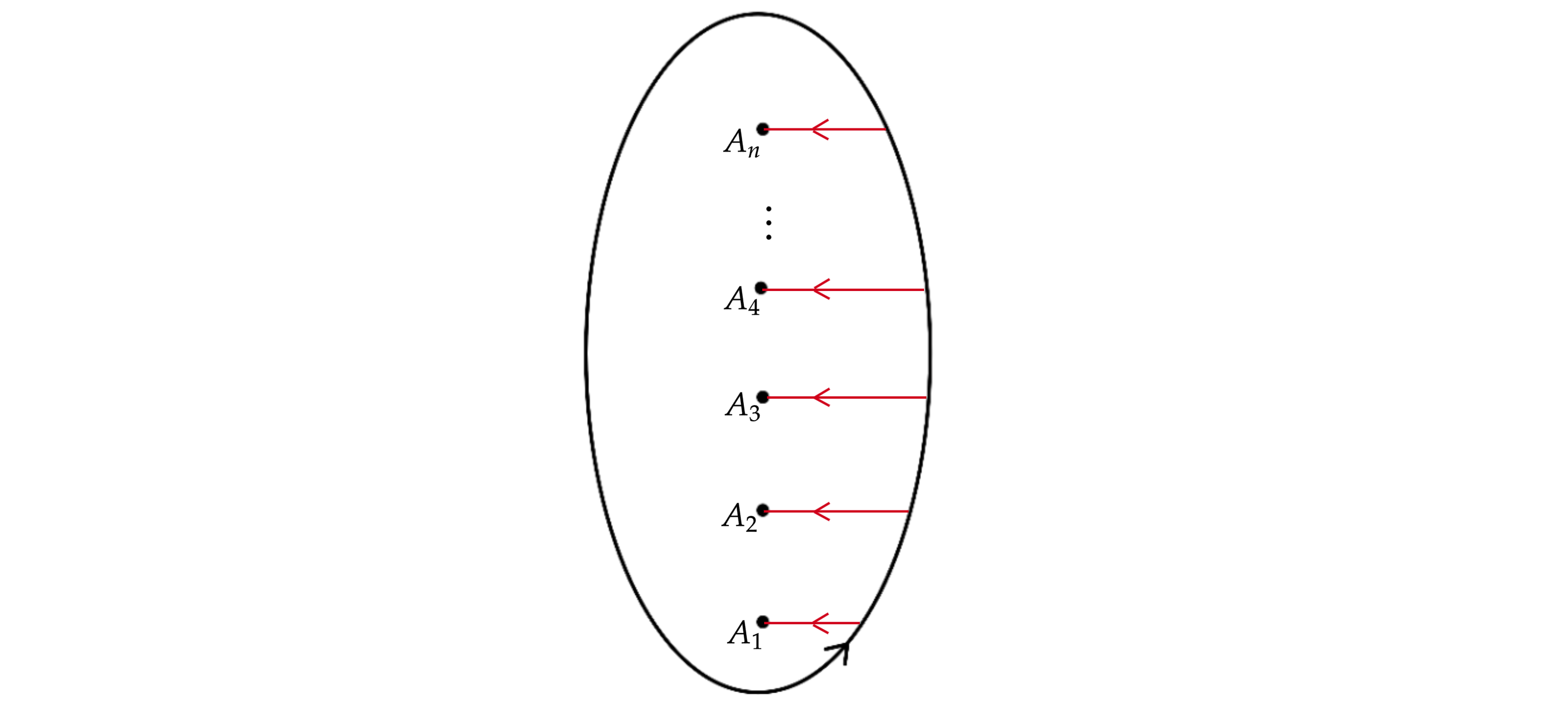}
    \caption{Branch cuts}
    \label{BranchCuts}
\end{figure}

\subsection{Contact branched coverings}\label{contact branched coverings}

Our main focus is the transverse-universality problem for transverse knots in the standard contact structure $(\mathbb{R}^3,\xi_{std})$, which is contactomorphic to $(\mathbb{R}^3,\xi_{rot})$. In the latter, any transverse knot can be represented as a braid, as previously discussed. This allows us to compute the self-linking number invariant using the formula $sl(L)=n_+-n_--b$, where $n_+$ and $n_-$ denote the number of positive and negative crossings, respectively, and $b$ is the number of strands \cite{Ver}. Recall that we can always decrease the self-linking number through a process called \textit{stabilization}; the converse is called \textit{destabilization}, and cannot always be done.

Related to this, we have the following fundamental result:

\begin{theorem}[{\cite[Theorem~4.2.8]{Cas}}]\label{MaxSl}
    Given any stabilized transverse knot, every contact branched covering with branch locus the given knot will be overtwisted.
\end{theorem}

If there is only one transverse representative for each possible self-linking number (up to the maximal one) in a fixed topological type of a knot, we say that it is \textit{transversely simple}. 

Now, we will outline the elements needed to prove Theorem \ref{Maintheo}. A detailed description of this construction is contained in \cite[Chapter~4]{Cas}.

To describe branched coverings over transverse links in $(\mathbb{R}^3,\xi_{rot})$, take the open book decomposition $(D^2,id)$ on $\mathbb{S}^3$ and the induced one on $\mathbb{R}^3$ described as before. Note that if $L$ is a transverse link braided with respect to these pages and having braid index $n$, the corresponding braid word determines, up to isotopy, a self-diffeomorphism of the $n$-punctured disk, which extends to a diffeomorphism isotopic to the identity of the filled disk. Fix $\Sigma_0$ as the initial disk for the braid word, and let $A_1,\dots,A_n$ its intersection points with $L$, such that the positive half twist $\sigma_i$ interchanges $A_i$ and $A_{i+1}$. Now, define a \textit{branch cut} associated with $A_i$ as a segment in $D^2$ that connects $A_i$ and $\partial D^2$ without passing through other points, and denote this branch cut by the same symbol, $A_i$. We will also assume that $\partial D^2$ with a fixed orientation intersects positively with each branch cut in the same order as the associated points, see Figure \ref{BranchCuts}.

\begin{figure}
    \centering
    \includegraphics[width=12cm]{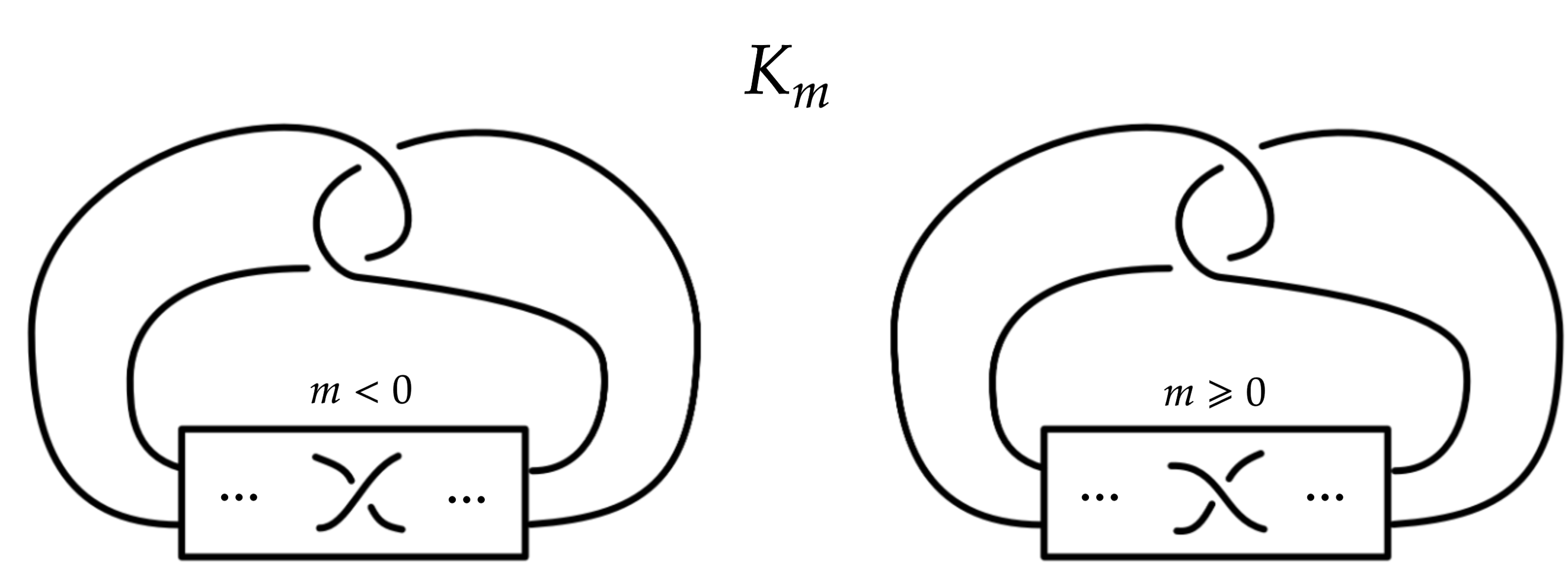}
    \caption{Twist knot $K_m$}
    \label{twist}
\end{figure} 

Let $\pi:(M,\xi)\rightarrow (\mathbb{S}^3, \xi_{std})$ be a contact branched covering, and consider again the open book decomposition $(D^2,\phi=id)$ for the base space. Then $D^2$ lifts to a surface $\tilde{\Sigma}$ through $\pi$, and we can consider $\tilde\phi$ as the gluing map in an associated open book decomposition of the total space, i.e., $\tilde\phi\in Aut(\tilde\Sigma)$ is such that $\tilde\phi=id$ on $\partial\tilde\Sigma$ and satisfies $\pi\circ\tilde\phi=id\circ\pi$ (see \cite[Section~3.4]{Cas}).
Take $\gamma$ an arc in $D^2$ that connects two points in the boundary, and encloses only a fixed $A_i$ and its branch cut clockwise, without intersecting any branching point, and let $\tilde\gamma$ be a lift of $\gamma$ with a fixed starting point. We have the following definition:

\textbf{Definition:} A \textit{branching word of $\phi(\gamma)$} is a word in the $A^{\pm 1}_1, \dots, A^{\pm 1}_n$ that describes the order in which this curve passes through the branch cuts, and the  sign of the exponent is the sign of the intersection. For the lift $\widetilde{\phi(\gamma)}$ we have the same construction for any fixed component. If the branching word is obtained from a minimal intersection with the branch cuts, we say that the branching word is reduced. 

\textbf{Remark:} Since far from the branching set $\pi$ is a covering map, the branching word for each component $\widetilde{\phi(\gamma)}$ is the same as for $\phi(\gamma)$.

 A \textit{detailed branching word} is a branching word for $\widetilde{\phi(\gamma)}$ where we add the information of the sheets. More precisely, if $x\in D^2$ is such that $x\in\gamma$ and $\pi^{-1}(x)=\{x_1,\dots,x_n\}$, where the subscripts represent an enumeration for the sheets, and  $\widetilde{\phi(\gamma)}_i$ is the unique lift of $\phi(\gamma)$ which starts at $x_i$, then we can trace the intersections of $\widetilde{\phi(\gamma)}_i$ with the branch cuts, for example, if it passes from the $s$-sheet to the $t$-sheet through the $A_j$ branch cut, the associated term is $A^{\pm 1}_{j,st}$ with the sign of the intersection. Note that a detailed branching word can be reduced and changed using isotopies of the curve as the following example shows: 
\begin{figure}
    \centering
    \includegraphics[width=14cm]{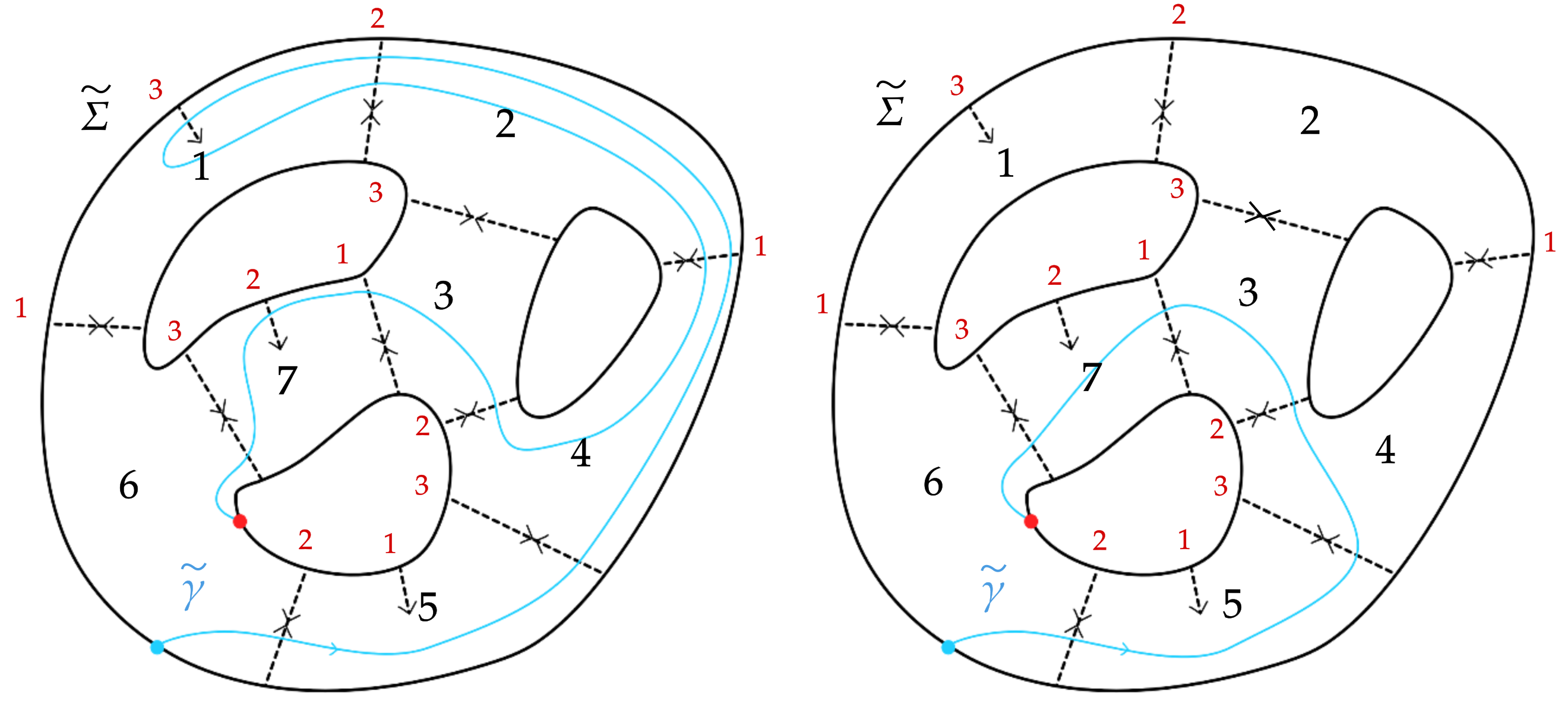}
    \caption{Isotopy of $\tilde\gamma$ in the surface $\tilde\Sigma$}
    \label{Example}
\end{figure}

\textbf{Example:} Consider the surface $\tilde\Sigma$ in Figure \ref{Example}. This surface can be viewed as a branched covering over $D^2$, where the numbers in black are a fixed enumeration for the sheets, and the ones in red represent the branch cuts, used to recover $\tilde\Sigma$. More precisely, this is a $7$-sheeted branched covering of $D^2$ with branch locus $\{A_1,A_2,A_3\}$ and monodromy:
\begin{align*}
        a_3&=(2\hspace{0.3cm}3)(4\hspace{0.3cm}5)(6\hspace{0.3cm}7)\\
        a_2&=(1\hspace{0.3cm}2)(3\hspace{0.3cm}4)(5\hspace{0.3cm}6)\\
        a_1&=(1\hspace{0.3cm}6)(2\hspace{0.3cm}4)(3\hspace{0.3cm}7).
    \end{align*}

where $a_i$ denotes the image in $S_7$ of the trivial generator of the fundamental group, corresponding to a loop based at the boundary that encloses only the branching point $A_i$ on its right side.

If we consider the curve $\tilde\gamma\subseteq\tilde\Sigma$ in the left-hand side of Figure \ref{Example}, which starts in the blue dot, we trace this curve using the following branching word:
$$
A_2A_3A_1A_2A_3A_2^{-1}A_1^{-1}A_2A_1A_2A_3^{-1}
$$
and its detailed branching word:
$$
A_{2,65}A_{3,54}A_{1,42}A_{2,21}A_{3,11}A_{2,12}^{-1}A_{1,24}^{-1}A_{2,43}A_{1,37}A_{2,77}A_{3,76}^{-1}
$$

Now, before we simplify this detailed branching word, we will mention the set of rules needed to do this process in general:
\begin{figure}
    \centering
\includegraphics[width=1\linewidth]{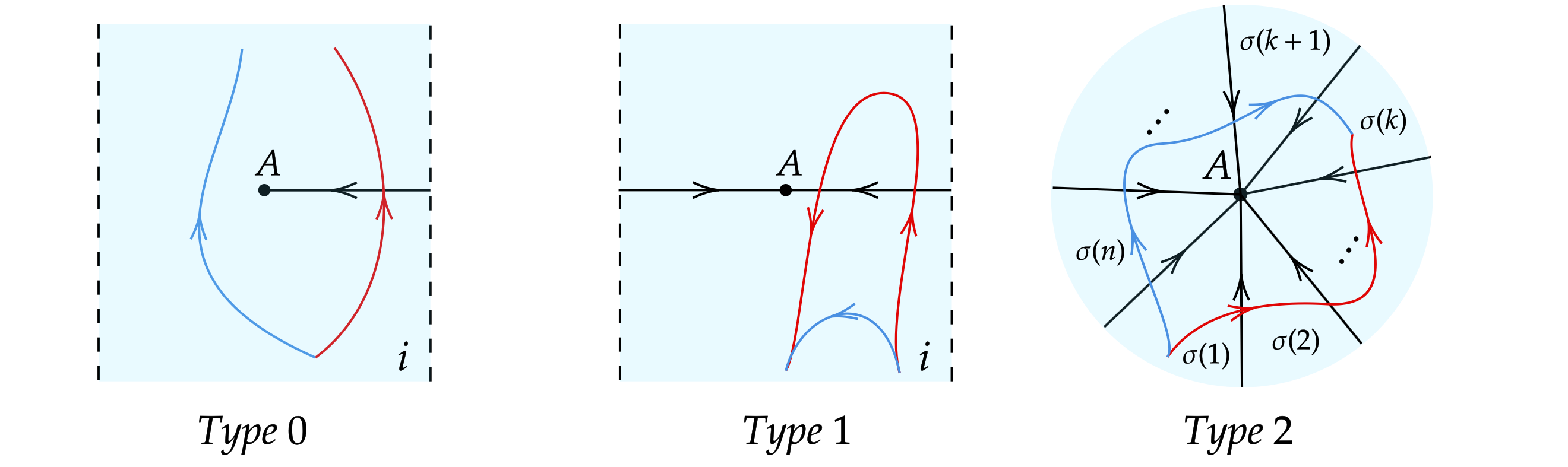}
    \caption{Different types of simplifications around a branching point $A$ with $s=1$. The red curve is the original, and the blue one is the one obtained after the isotopy. The enumerations correspond to the specific sheet}
    \label{Simplifications}
\end{figure}
there are different types of simplifications that can be done in a detailed branching word to reduce it. Note that any isotopy of a curve inside any sheet does not change the detailed branching word, even more, any simplification have to come from the behavior of the curve near a branching point. We will separate them into 3 types:

Fix a branching point $A$ and $s\in\{-1,1\}$:
\begin{itemize}
    \item \textbf{Type 0:}  
    $$
A_{ii}^{s}\Rightarrow \emptyset
    $$
    \item \textbf{Type 1:}
    $$
A_{ij}^sA_{ji}^{-s}\Rightarrow \emptyset
    $$
    \item \textbf{Type 2:} If $A$ is a branching point of index $n$ with associated permutation $$(\sigma(1)\hspace{0.3cm}\sigma(2)\hspace{0.3cm}\dots\hspace{0.3cm}\sigma(n))$$
    then, for $k>\lceil{\frac{n}{2}}\rceil$:
    $$
A_{\sigma(1)\sigma(2)}^{s}A_{\sigma(2)\sigma(3)}^{s}\dots A_{\sigma(k-1)\sigma(k)}^{s}\Rightarrow A_{\sigma(1)\sigma(n)}^{-s}A_{\sigma(n)\sigma(n-1)}^{-s}\dots A_{\sigma(k+1)\sigma(k)}^{-s}
    $$
\end{itemize}
Figure \ref{Simplifications} shows these simplifications for $s=1$.

\textbf{Remark:}
    The type $0$ simplification is a particular case of type $2$ with $n=k=1$, but since type $0$ and $1$ are going to be the most common ones in our future proofs we separated it from the other.

Using this in the previous example, we can simplify the detailed branching word:

\begin{align*}
A_{2,65}A_{3,54}A_{1,42}A_{2,21}A_{3,11}&A_{2,12}^{-1}A_{1,24}^{-1}A_{2,43}A_{1,37}A_{2,77}A_{3,76}^{-1}\\
   A_{2,65}A_{3,54}A_{1,42}A_{2,21}A_{2,12}^{-1}&A_{1,24}^{-1}A_{2,43}A_{1,37}A_{3,76}^{-1}\\
   A_{2,65}A_{3,54}A_{1,42}&A_{1,24}^{-1}A_{2,43}A_{1,37}A_{3,76}^{-1}\\
   A_{2,65}A_{3,54}A_{2,43}&A_{1,37}A_{3,76}^{-1}
\end{align*}

Which is equivalent to the detailed branching word of the right-hand side of Figure \ref{Example}.

 The branching words are used to determine if the monodromy map is right- or left-veering. For this, if we consider two non-isotopic curves, $\gamma_1,\gamma_2:[0,1]\rightarrow \Sigma$, on a compact, oriented surface $\Sigma$ with non-empty boundary, with the same initial point on $\partial \Sigma$, we say that $\gamma_1$ is \textit{to the right} of $\gamma_2$ if, after isotoping them rel boundary for minimal intersection, the induced orientation from the pair $\{\dot\gamma_1(0),\dot\gamma_2(0)\}$ agrees with the one in $\Sigma$. If the induced orientation is the opposite, we say that $\gamma_1$ is \textit{to the left} of $\gamma_2$.
 
  Recall that an automorphism of a surface $\phi$ that is the identity on the boundary is \textit{right-veering} if, for every properly embedded arc $\gamma$ in the surface, $\phi(\gamma)$ is to the right or is isotopic to $\gamma$. We call it \textit{left-veering} if it is not right-veering. Using this terminology and with the same notation, the following lemma can be proved:

\begin{lemma}[{\cite[Lemma~4.2.7]{Cas}}]
    If $\gamma$ is a properly embedded arc that encloses only the branch cut $A_i$ at its right side, such that for some connected component $\alpha$ of $\tilde\gamma$, the reduced branching word of $\tilde\phi(\alpha)$ does not start with $A_i$ and is not empty, then $\tilde\phi$ is left-veering.
\end{lemma}

\subsection{Twist knots}

Let $m\in\mathbb{Z}$. The twist knot $K_m$ can be obtained by linking the left and right parts of a trivial knot and performing $m$ half-twists, which could be positive or negative, according to the sign of $m$; this is represented in Figure \ref{twist}. An important fact about some members of this family is the following:

 \begin{theorem}[{\cite[Section~7.1]{Pu}}] $K_m$ is hyperbolic if and only if $m\geq2$ or $m\leq -3$.
 \end{theorem}

 It follows that for those values of $m$, $K_m$ is non-toroidal, and since every nontrivial $K_m$ is a 2-bridge knot, $K_m$  will be universal for $m\geq2$ and $m\leq-3$ \cite{HiLoMo3}.

Related to transverse representatives of  twist knots, we have the following result:

 \begin{theorem}[\textbf{Classification of transverse twist knots} \cite{EtLenVer}]\label{classification}
\begin{itemize}\
    \item If $m$ is odd, then $K_m$ is transversely simple. Moreover, the transverse representative of $K_m$ with maximal self-linking number has $sl=-m-4$ if $m>-1$, $sl=-1$ if $m=-1$ and $sl=-3$ if $m<-1$.
    \item If $m\geq-2$ is even, then $K_m$ is transversely simple and every transverse representative with maximal self-linking number has $sl=-m-1$.
    \item If $m\leq-4$ is even, then $K_m$ is transversely non-simple or $m=-4$. If $m=-2n$, there are $\lceil \frac{n}{2} \rceil$ distinct transverse representatives of $K_m$ with maximal self-linking number $sl=1$.
\end{itemize}
 \end{theorem}

\begin{figure}
    \centering
    \includegraphics[width=14cm]{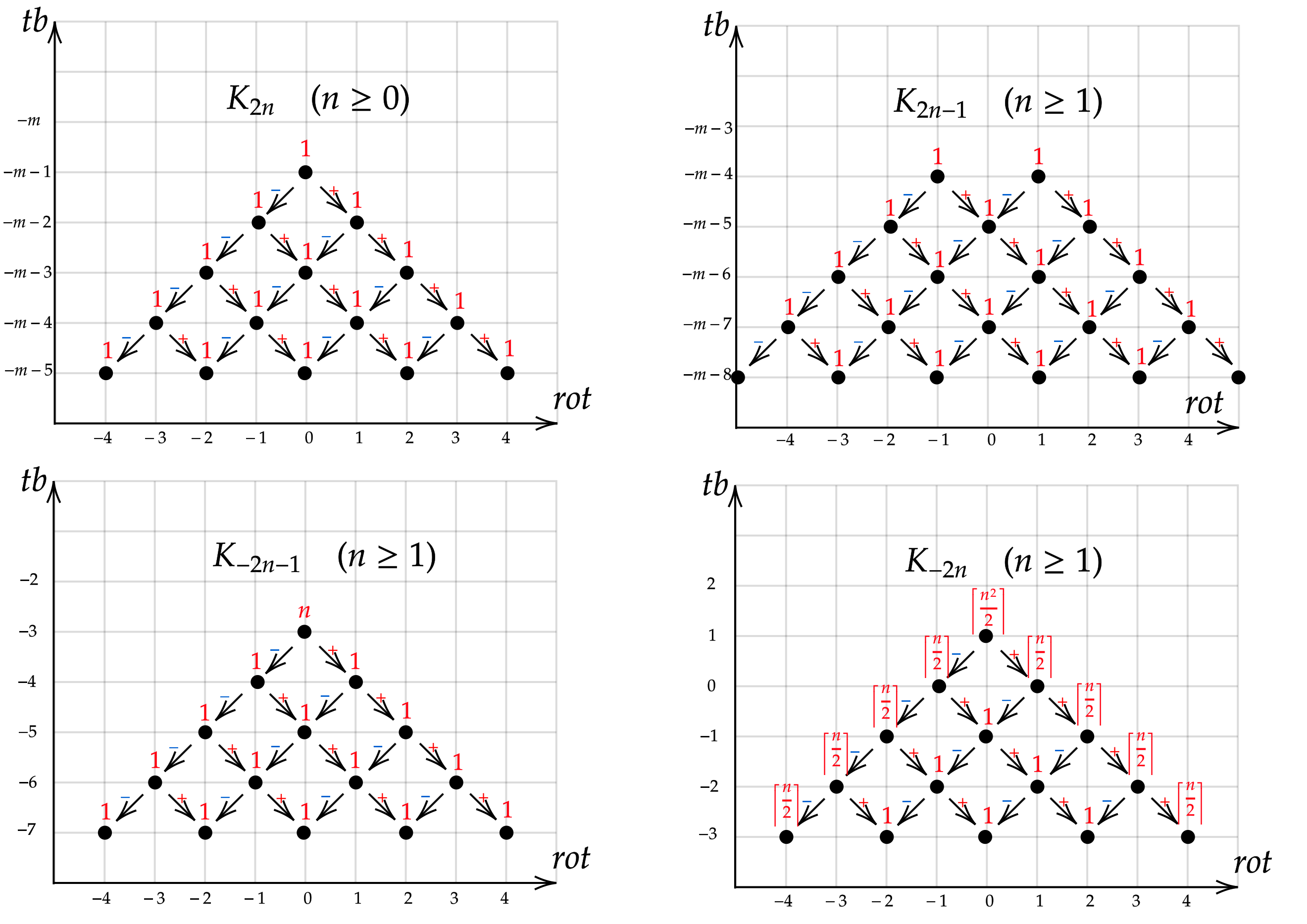}
    \caption{Legendrian mountain range diagrams for $K_m$, $m\neq -1$. The numbers indicate de number of non-equivalent Legendrian representatives for the given pair $(rot,tb)$.}
    \label{Legendrian mountain range for twist}
\end{figure}

Recall that any transverse link can be obtained as the positive transverse push-off of a Legendrian link \cite{EpsFuMe}. And consider Figure~\ref{Legendrian mountain range for twist}, which depicts the Legendrian mountain range for Legendrian twist knots (Taken from \cite{EtLenVer}). In this context, the following result holds:

\begin{theorem}[{\cite[Theorem~3.3]{EtLenVer}}]
    Let $K$ be a Legendrian representative of the twist knot $K_m$. If $K$ does not realize the maximal Thurston-Bennequin number among all Legendrian representatives of $K_m$, then $K$ destabilizes. 
\end{theorem}

Observe that Figure~\ref{Legendrian mountain range for twist} contains only global maxima, and each vertex connected by four arrows corresponds to a unique Legendrian class. It follows that the non-destabilizable transverse twist knots are precisely the positive transverse push-offs of the Legendrian knots located at the left-most part of each diagram. Equivalently, these knots realize the maximal self-linking number among all possible positive transverse push-offs. Hence, by Theorem~\ref{MaxSl}, we may restrict our attention throughout the proof to transverse representatives with maximal self-linking number.
 
\section{Proof of the main theorem}

We now proceed with the proof of Theorem~\ref{Maintheo}, distinguishing three cases according to the value of $m$. Although the arguments share a similar structure, the cases are not equivalent. Let 
$\pi : (M,\xi) \longrightarrow (\mathbb{S}^3,\xi_{std})$
be a non-trivial contact branched covering of degree $k_0$, branched along a transverse twist knot $K_m$ with maximal self-linking number.

\textbf{Case 1: $m\geq2$ even:}

\begin{figure}[h]
    \centering
    \includegraphics[width=14cm]{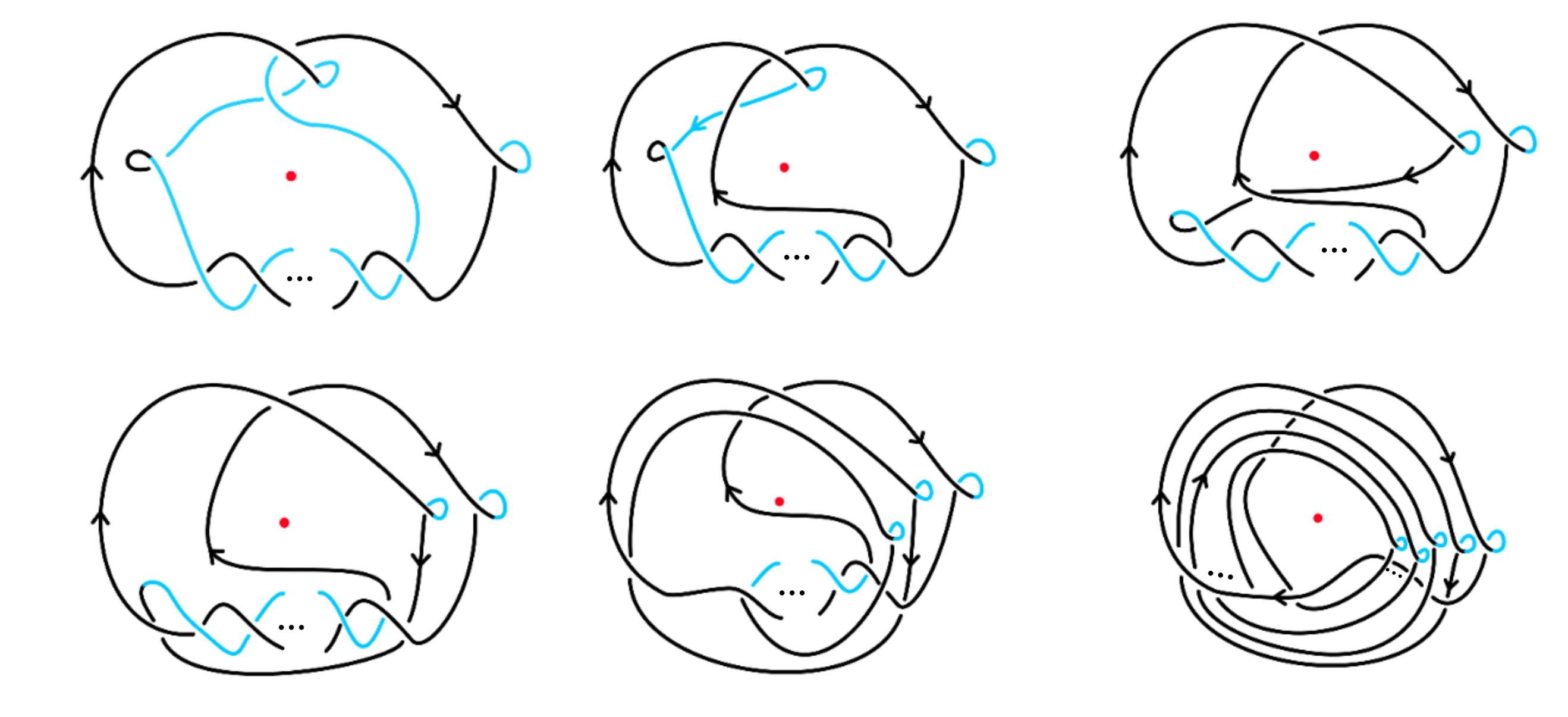}
    \caption{Braiding, case 1}
    \label{braiding1}
\end{figure}

By Theorem \ref{classification}, $sl(K_m)=-m-1$, and then, by transverse simplicity, we only need to take a transverse knot in the same topological class of $K_m$ with the given self-linking number. Consider the given transverse knot and the braiding process in Figure \ref{braiding1}, where the leftmost top knot is a front projection of $K_m$. From left to right and top to bottom, we present the transverse isotopies that transform $K_m$ into its braid form.

\begin{figure}
    \centering
    \includegraphics[width=14cm]{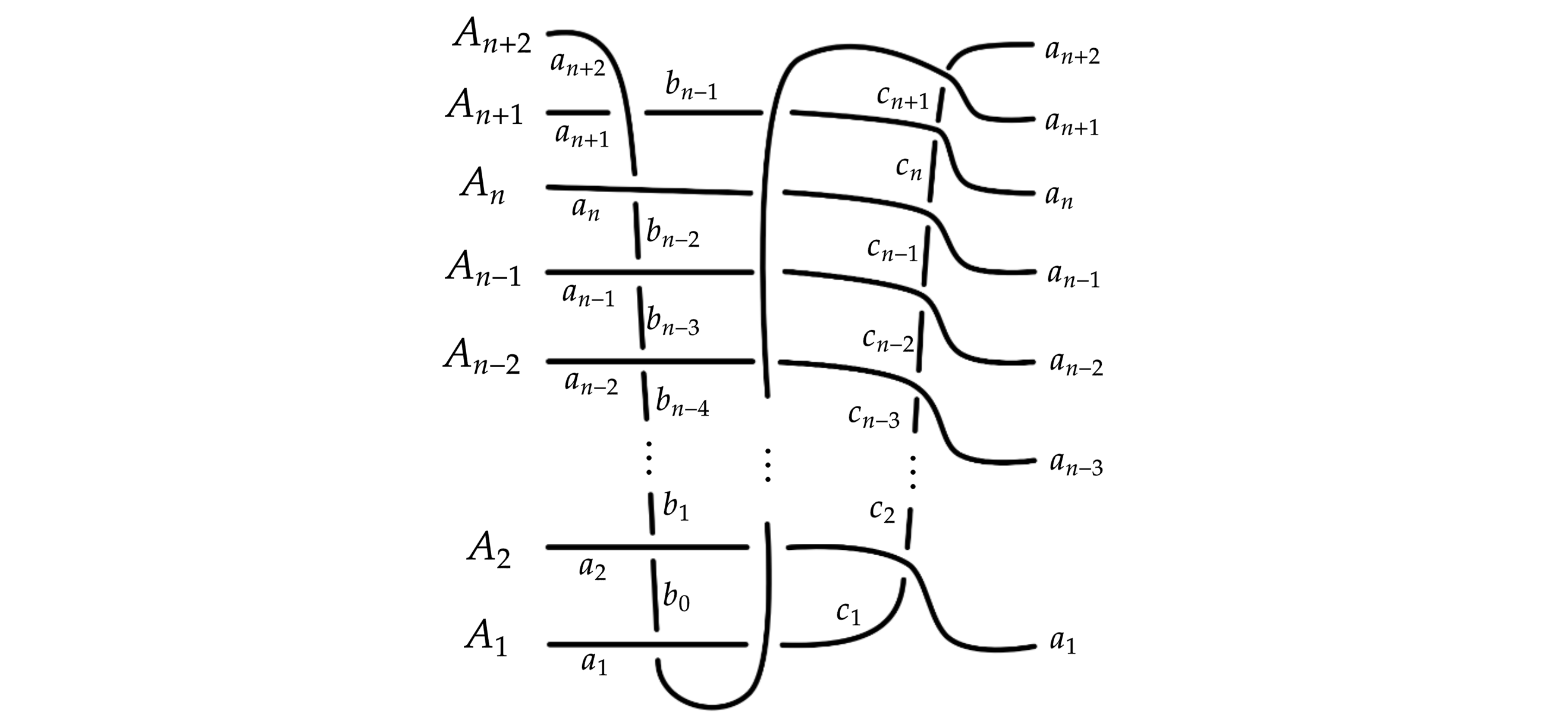}
    \caption{Braid coloration, case 1}
    \label{coloration1}

\end{figure}

Let $n=\frac{m}{2}$. Figure \ref{coloration1} shows the induced braid from the above process and fixes a strand enumeration $A_1,A_2,\dots, A_{n+2}$, such that the braid word has the following expression:
$$
\sigma_{n+1}\sigma_{n}^{-1}\sigma_{n-1}^{-1}\dots\sigma_{2}^{-1}\sigma_{1}^{-2}\sigma_{2}^{-1}\dots\sigma_{n+1}^{-1}\sigma_{1}\sigma_{2}\dots\sigma_{n+1}
$$

\begin{figure}
    \centering
\includegraphics[width=15cm]{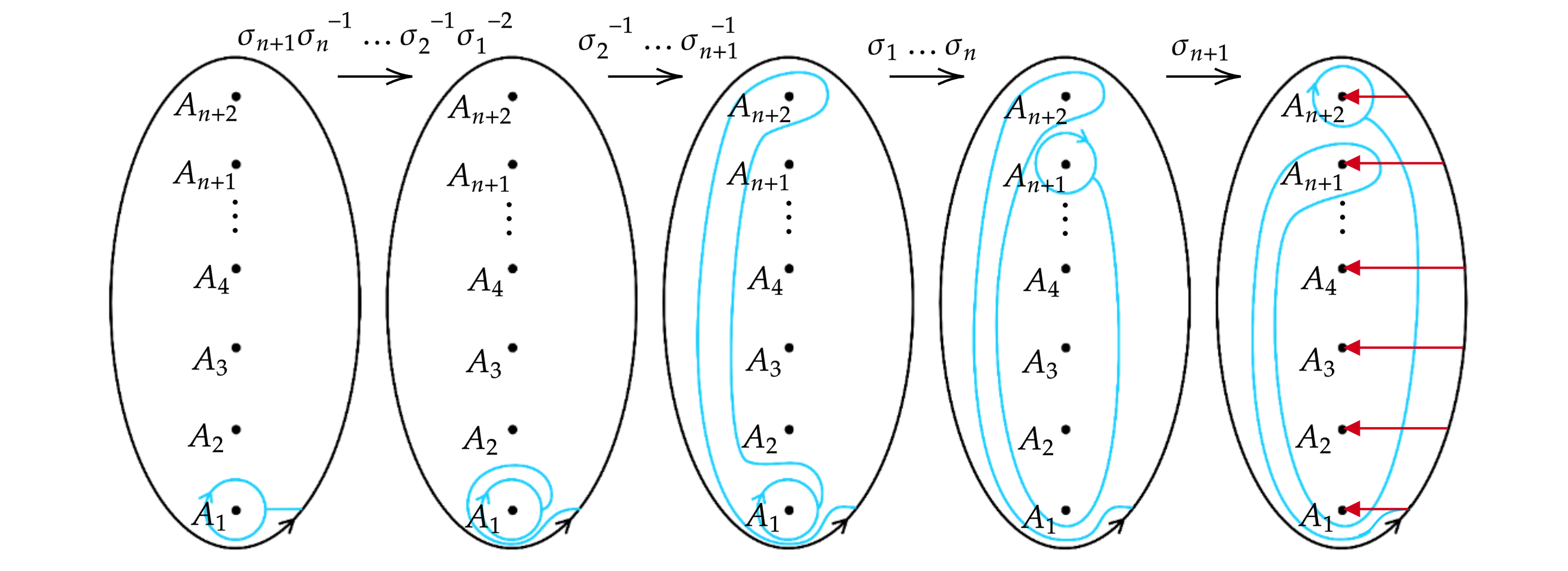}
    \caption{Deduction of the branching word, case 1}
    \label{branching1}
\end{figure}

Using this, we can recover the branching word associated with the $A_1$ strand, as Figure \ref{branching1} shows. It follows that it takes the form $$W=A_{n+1}^{-1}A_1A_2\dots A_{n+1}A_{n+2}^{-1}A_{n+1}^{-1}\dots A_{1}^{-1}A_{n+1}A_1$$

To simplify the notation, let $\Omega=A_1A_2\dots A_{n+1}$, then $W=A_{n+1}^{-1}\Omega A_{n+2}^{-1}\Omega^{-1}A_{n+1}A_1$.

Recall that $\Sigma_0$ is the initial page of the open book decomposition of $\mathbb{S}^3$, just as in Section \ref{contact branched coverings}. Denote with $a_j$ the coloring of $A_j\in\Sigma_0$; that is, $a_j$ is the image under the monodromy map (induced by the branched cover) of the associated class of a meridional curve for the respective strand, see Figure \ref{coloration1}. 
 
Let $\gamma \subseteq D^2$ denote the blue curve in the first disk of Figure~\ref{branching1}. Observe that $\gamma$ is a properly embedded arc enclosing the branch cut associated with $A_1$ on its right side. To verify that $\tilde{\phi}$ is left-veering, as mentioned previously, it suffices to show that there exists a reduced detailed branching word associated with $\tilde{\gamma}$ that does not begin with $A_1$ or is empty. Equivalently, we must prove that the branching word cannot be reduced either to another word beginning with $A_1$ or to the empty word. For the specific case of $W$, we need to verify that the initial term $A_{n+1}^{-1}$ and the almost-middle $A_{n+2}^{-1}$ cannot be eliminated.  This is equivalent to the existence of an $i$ for which:
 
\begin{equation}\tag{C1}\label{cond1}
\begin{matrix} a_{n+1}(i)\neq i\\
    a_{n+2}^{-1}a_{n+1}\dots a_1(i)\neq a_{n+1}\dots a_{1}(i)
\end{matrix}
\end{equation}
 
     Since $\pi$ is a covering map at the complement of the branching points, then $\pi^{-1}(\gamma)$ has $k_0$ components. Let $\alpha_i$ be the component of $\tilde\gamma:=\pi^{-1}(\gamma)$ that has its starting point in the $a_{n+1}(i)$-sheet. We need to study $\tilde\phi(\alpha_i)$ such that $i$ satisfies the conditions \ref{cond1}.

\textbf{Claim:} There exists an $i$ that satisfies condition  \ref{cond1}.

By definition, $W$ is the branching word of $\tilde\phi(\alpha_i)$ for all $i$. To prove this claim, recall first that the initial condition in \ref{cond1} tells us that $\tilde\phi(\alpha_i)$ branches when it passes the branch cut associated to $A_{n+1}$; therefore, since $A_{n+1}$ is a branching point, there exists at least one $i$ with $a_{n+1}(i)\neq i$. Suppose now by contradiction that for all of those $i$'s, the second condition in \ref{cond1} is false. If we consider the detailed branching word associated with $\tilde\phi(\alpha_i)$ with $j=a_{n+1}(i)$, it has the following form:
$$
A_{n+1,ji}^{-1}\Omega_{ik}A_{n+2,ks}^{-1}\Omega^{-1}_{st}A_{n+1,tu}A_{1,uj}
$$
The first condition implies that $i\neq j$, and since we are assuming that the second condition is false, $k=s$, therefore we can reduce this detailed branching word:
\begin{align*}
A_{n+1,ji}^{-1}\Omega_{is}A_{n+2,ss}^{-1}&\Omega^{-1}_{st}A_{n+1,tu}A_{1,uj}\\
A_{n+1,ji}^{-1}\Omega_{is}&\Omega^{-1}_{st}A_{n+1,tu}A_{1,uj}\\
A_{n+1,ji}^{-1}&A_{n+1,iu}A_{1,uj}\\
&A_{1,jj}\\
&\emptyset
\end{align*}
It follows that if $A_{n+1}$ branches the $i$-sheet, so that $a_{n+1}(i)=j$ (i.e., the branch cut $A_{n+1}$ connects the $i$-sheet with the $j$-sheet), and no branching occurs at $A_{n+2}$ when passing through it, then $A_1$ does not branch the $j$-sheet. Since $j=a_{n+1}(i)$ is an arbitrary nonfixed element of $a_{n+1}$, $a_{1}$ and $a_{n+1}$ are two disjoint permutations. For that reason, those two permutations commute.

Consider now the braid coloration in Figure \ref{coloration1}. By Wirtinger's relations for knots, we have the following relations: 
     \begin{align*}
        \tag{1} b_{n-2}&=a_na_{n+2}a_{n}^{-1} \\
        \tag{2}b_i&=a_{i+2}b_{i+1}a_{i+2}^{-1} &0\leq i\leq n-3\\
        \tag{3}a_{n+1}&=a_1b_0a_1^{-1}\\
        \tag{4}a_i&=a_{n+1}a_{i+1}a_{n+1}^{-1} &1\leq i\leq n-1\\
        \tag{5}c_1&=a_{n+1}a_1a_{n+1}^{-1}\\
        \tag{6}a_{n+2}&=a_{n+1}c_{n+1}a_{n+1}^{-1}\\
        \tag{7}c_i&=a_{i-1}^{-1}c_{i-1}a_{i-1} &2\leq i \leq n+1
    \end{align*}

We already know that $a_{n+1}a_1=a_1a_{n+1}$. With this and the previous relations, it is easy to verify that the image of $\pi_1(\mathbb{S}^2-K_m)$ under the monodromy map is generated by $a_1$ and $a_{n+1}$. Even more:
\begin{align*}
    a_i&=a_1\hspace{1cm} 1\leq i\leq n\\
    b_i&=a_{n+1} \hspace{0.6cm} 0\leq i\leq n-2\\
    c_i&=a_1\hspace{1cm} 1\leq i\leq n+1\\
    a_{n+2}&=a_1
\end{align*}
This, together with the relation $(1)$, imply that $a_{n+1}=a_1$,
but this is not possible since $a_{n+1}$ and $a_1$ are nontrivial permutations without nonfixed numbers in common. It follows that there must exist an $i$ that satisfies \ref{cond1}.

Consider again the detailed branching word associated to $\tilde\phi(\alpha_i)$, and fix an $i$ that satisfies \ref{cond1}:
$$
A_{n+1,ji}^{-1}\Omega_{ik}A_{n+2,ks}^{-1}\Omega^{-1}_{st}A_{n+1,tu}A_{1,uj}
$$
Since $j=a_{n+1}(i)\neq i$, the first term $A^{-1}_{n+1,ji}$ cannot be eliminated unless we cancel it with another one of the type $A_{n+1}^{\pm1}$. But, by \ref{cond1}:
$$a_{n+2}^{-1}a_{n+1}\dots a_1(i)\neq a_{n+1}\dots a_{1}(i)\Rightarrow k\neq s$$
this means that we cannot cancel $A^{-1}_{n+2,ks}$ and for that reason, $A^{-1}_{n+1,ji}$ could only be eliminated with the $A_{n+1}$ contained in $\Omega_{ik}$, but this can only happen if the intermediate terms are erased. In any case, the branching word does not start with $A_1$, and by this, we conclude that $\tilde\phi$ is left-veering; therefore, any contact branched cover over this knot is overtwisted.

\textbf{Case 2: $m\leq-3$ odd}

\begin{figure}
\centering
\includegraphics[width=16cm]{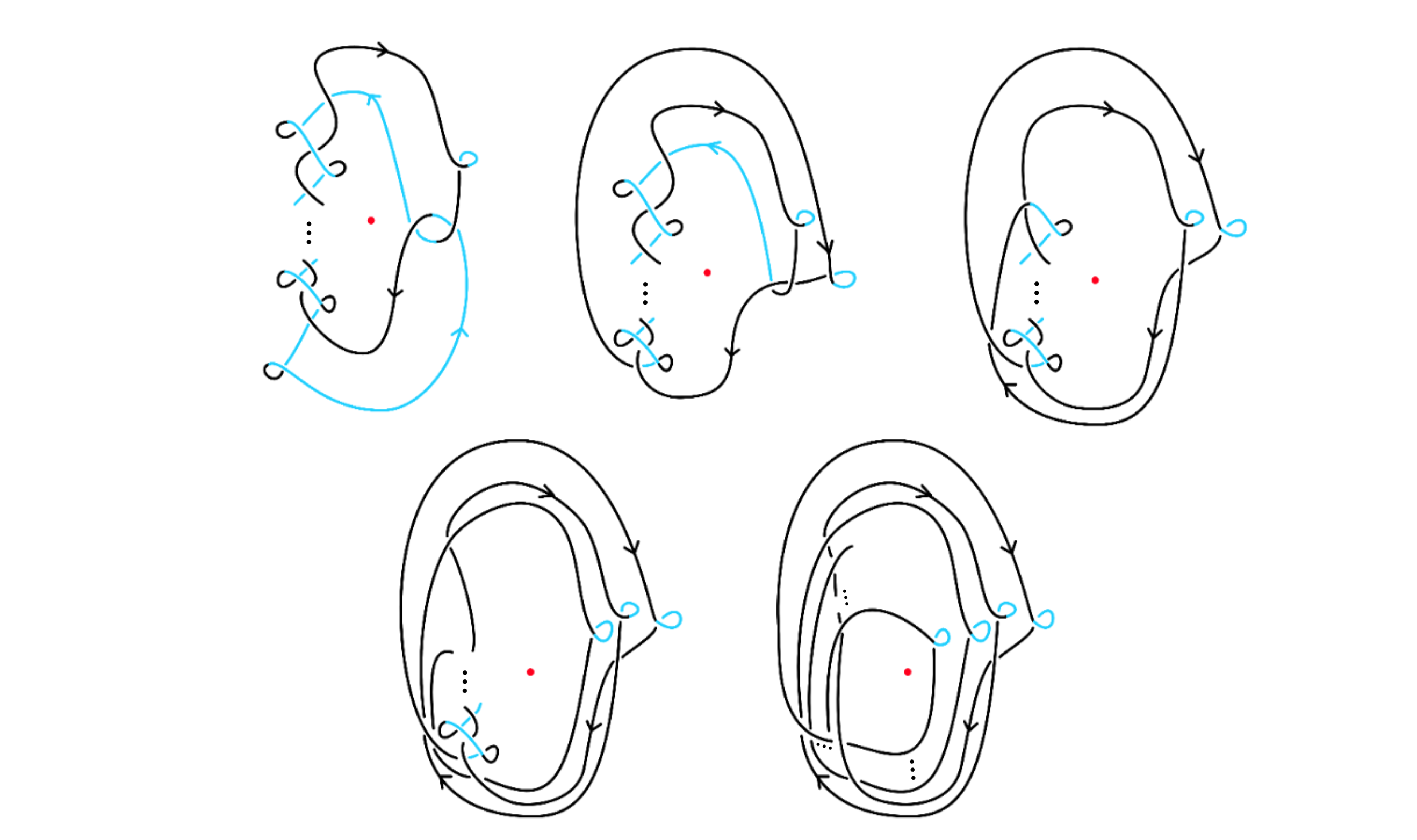}
\caption{Braiding, case 2}
\label{braiding2}
\end{figure}
\begin{figure}
         \centering
         \includegraphics[width=16cm]{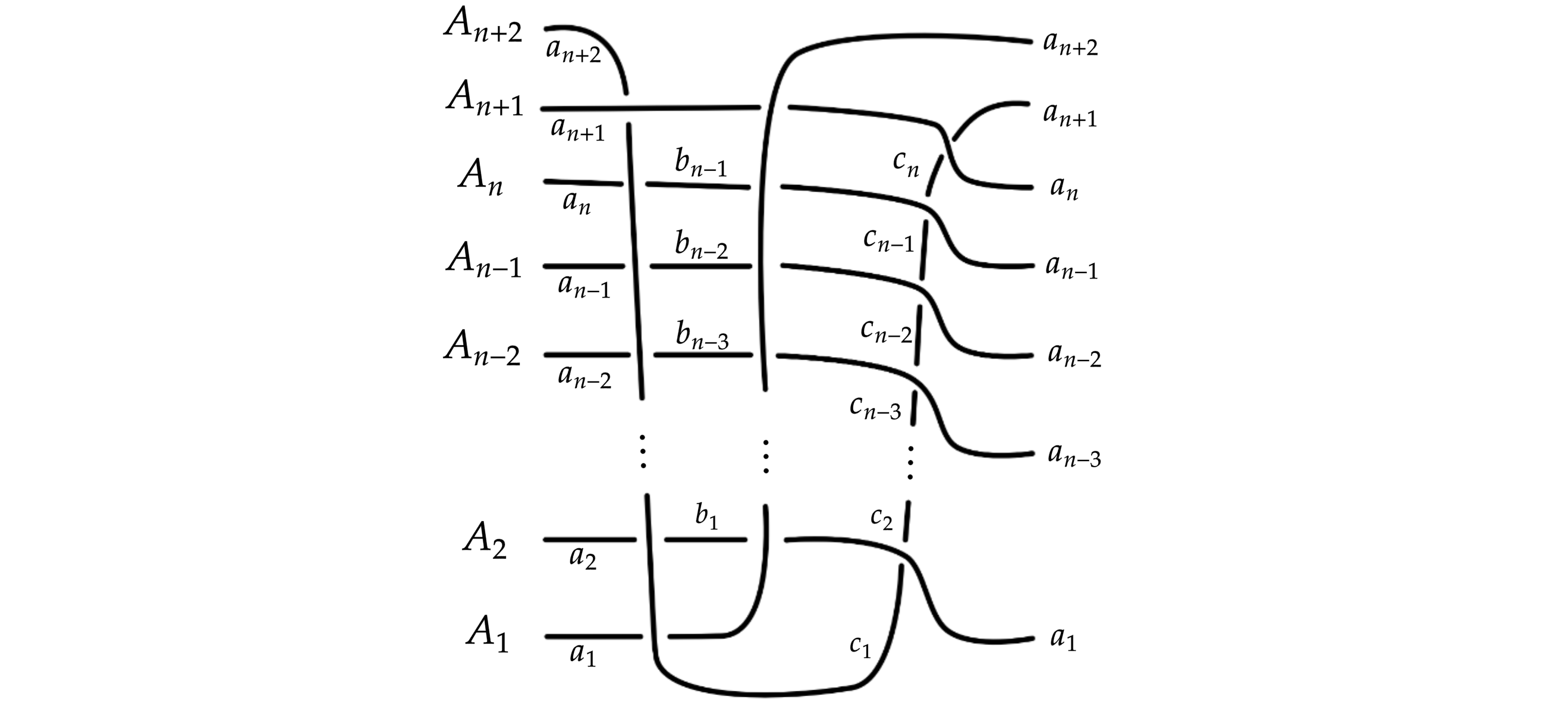}
         \caption{Braid coloration, case 2}
    \label{coloration2}
     \end{figure}

 
Let  $n=\frac{m-1}{2}$. For the transverse knot with $sl=-3$ and the given braiding process in Figure \ref{braiding2}, we have the braid word (see Figure \ref{coloration2}) :
$$
\sigma_{n+1}^{-1}\sigma_{n}\dots \sigma_{2}\sigma_1\sigma_{2}^{-1}\dots \sigma_{n+1}^{-1}\sigma_1\dots \sigma_{n}
$$

\begin{figure}
\centering
\includegraphics[width=16cm]{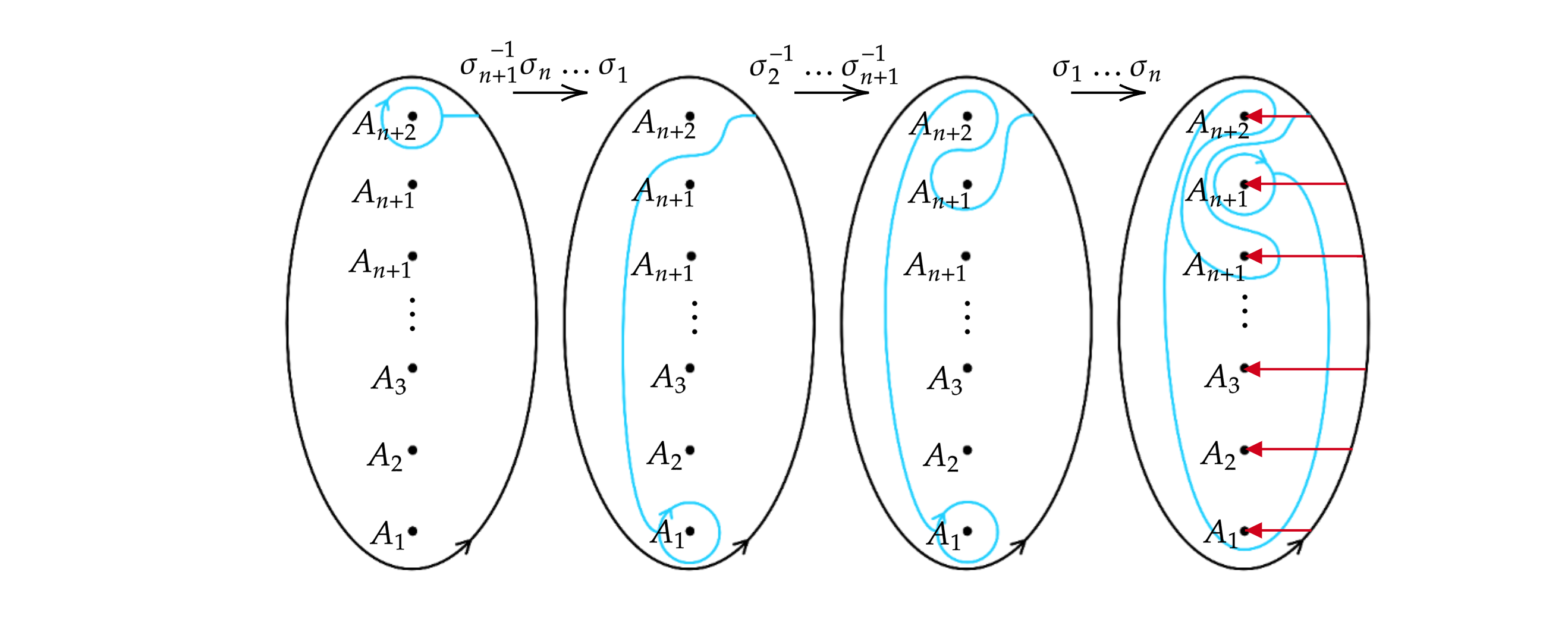}
\caption{Deduction of the branching word, case 2}
\label{branching2}
\end{figure}

And, calculating the branching word induced by the $A_{n+2}$-strand, as in Figure \ref{branching2}, we get $W=A_{n}^{-1}A_{n+2}A_{1}\dots A_{n}A_{n+1}^{-1}A_{n}^{-1}\dots A_{1}^{-1}A_{n+2}^{-1}A_{n}A_{n+2}$.

Defining again $\Omega=A_{n+2}A_{1}\dots A_{n}$, we have $W=A_{n}^{-1}\Omega A_{n+1}^{-1}\Omega^{-1}A_{n}A_{n+2}$.

\textbf{Claim:} There exists an $i$ such that:
\begin{equation}\tag{C2}\label{cond2}
\begin{matrix}
a_{n}(i)\neq i\\
a_{n+1}^{-1}a_{n}\dots a_1a_{n+2}(i)\neq a_{n}\dots a_{1}a_{n+2}(i)
\end{matrix}
\end{equation}

 As before, by definition, there exists at least one $i$ that satisfies the first part of the above condition, so suppose by contradiction that for each of those $i$, the second condition is actually an equality. In the detailed branching word
 $$
A_{n,ji}^{-1}\Omega_{ik} A_{n+1,ks}^{-1}\Omega_{st}^{-1}A_{n,tu}A_{n+2,uj}
 $$
 this is equivalent to say that $k=s$, and hence, after the reduction process, we obtain:
\begin{align*}
A_{n,ji}^{-1}\Omega_{is} A_{n+1,ss}^{-1}&\Omega_{st}^{-1}A_{n,tu}A_{n+2,uj}\\
A_{n,ji}^{-1}\Omega_{is}&\Omega_{si}^{-1}A_{n,iu}A_{n+2,uj}\\
A_{n,ji}^{-1}&A_{n,iu}A_{n+2,uj}\\
&A_{n+2,jj}\\
&\emptyset
\end{align*}

In particular, if $i\in a_{n+1}$, then $a_{n+1}(i)\notin a_{n+2}$, hence $a_{n+1}a_{n+2}=a_{n+2}a_{n+1}$. We again use the braid coloration of Figure \ref{coloration2} to obtain the relations:
     \begin{align*}
    \tag{1}c_1&=a_{n+1}a_{n+2}a_{n+1}^{-1}\\
    \tag{2}b_i&=c_1^{-1}a_{i+1}c_1& 1\leq i\leq n-1\\
    \tag{3} a_{n+2}&=c_1^{-1}a_1c_1\\
    \tag{4} a_i&=a_{n+2}b_i a_{n+2}^{-1} &1\leq i\leq n-1\\
    \tag{5} a_n&=a_{n+2}a_{n+1}a_{n+2}^{-1}\\
    \tag{6} c_i&=a_{i-1}^{-1}c_{i-1}a_{i-1} &2\leq i\leq n\\
    \tag{7} a_{n+1}&=a_{n}^{-1}c_na_n
    \end{align*}

Since $a_{n+1}a_{n+2}=a_{n+2}a_{n+1}$, again:
\begin{align*}
    a_i&=a_{n+2}\hspace{1cm} 1\leq i\leq n\\
    c_i&=a_{n+2}\hspace{1cm} 1\leq i\leq n\\
\end{align*}
After using this in $(5)$, it follows that $a_{n+1}=a_{n+2}$, which is again a contradiction. Therefore, we cannot eliminate the first term of the branching word, and hence $\tilde\phi$ must be left-veering.\hspace{0.3cm}\qedsymbol{}

\textbf{Case 3: $m\geq2$ odd}
\begin{figure}
         \centering
         \includegraphics[width=14cm]{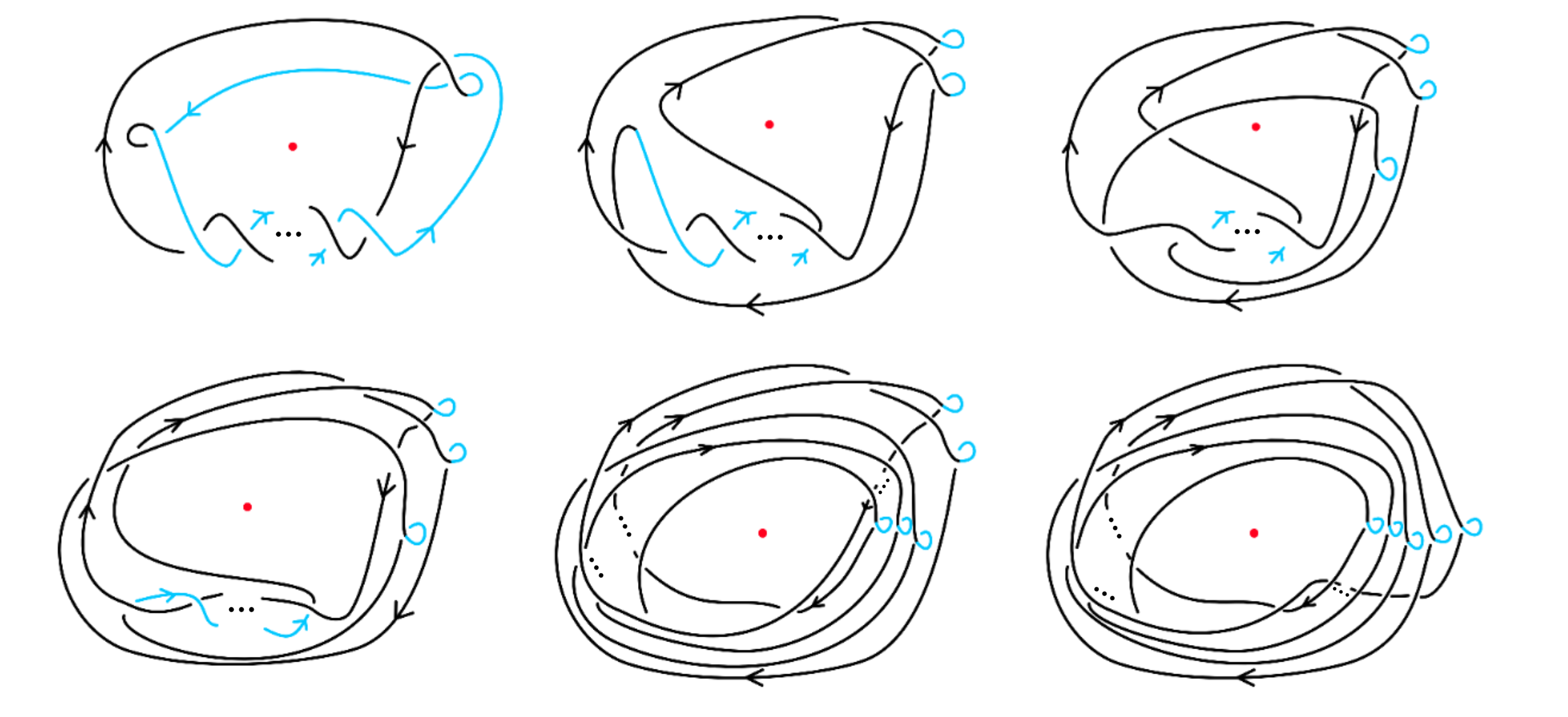}
         \caption{Braiding, case 3}
    \label{braiding3}
\end{figure}

\begin{figure}
         \centering
         \includegraphics[width=14cm]{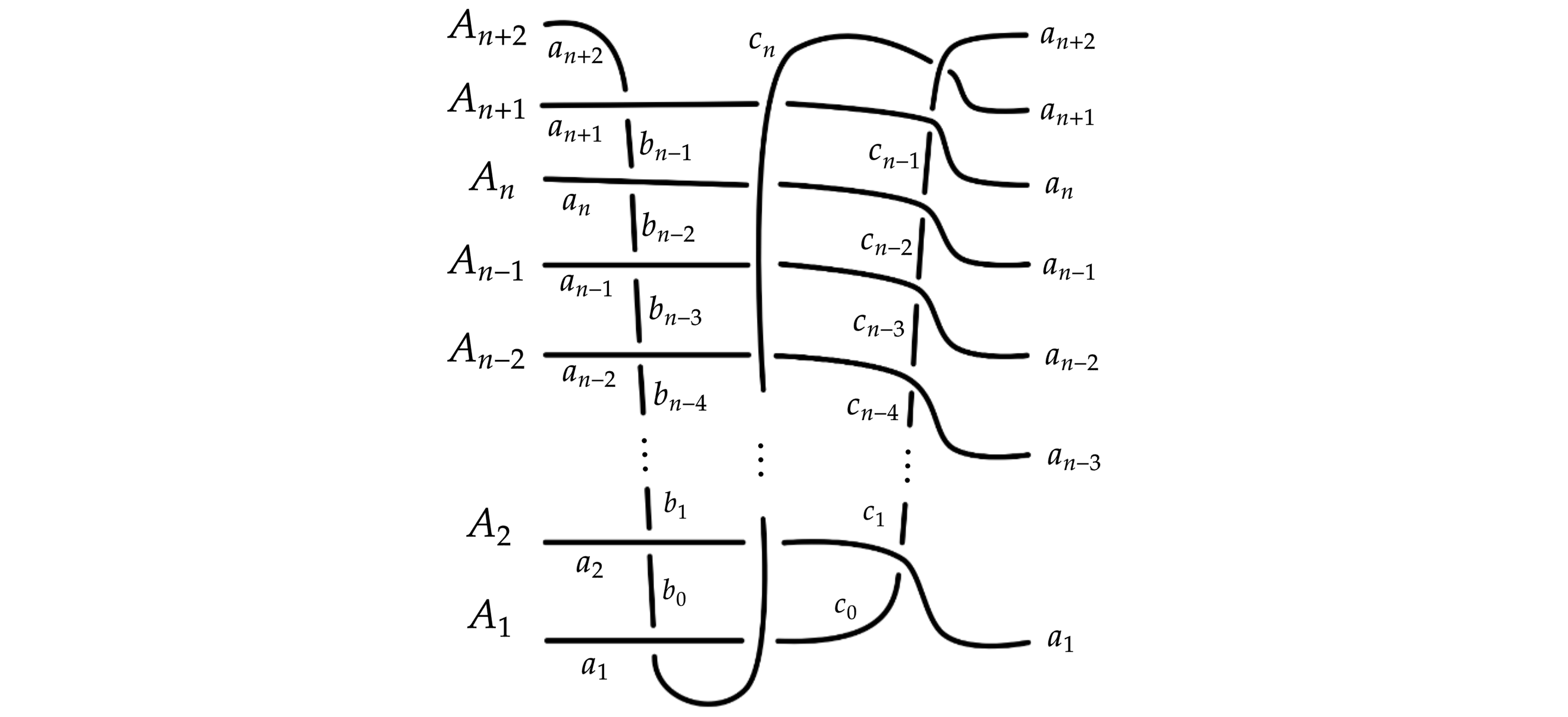}
         \caption{Braid coloration, case 3}
    \label{coloration3}
     \end{figure}

In this case, the transverse representative that we use has $sl=-m-4$. So, for example, take the one given in Figure \ref{braiding3}. Let $n=\frac{m-1}{2}$, as shown in Figure \ref{coloration3},
the braid word takes the form
$$
\sigma_{n+1}^{-1}\dots \sigma_2^{-1}\sigma_{1}^{-2}\sigma_{2}^{-1}\dots \sigma_{n+1}^{-1}\sigma_1\sigma_2\dots\sigma_{n}\sigma_{n+1}^{-1}
$$

 \begin{figure}
         \centering
         \includegraphics[width=16cm]{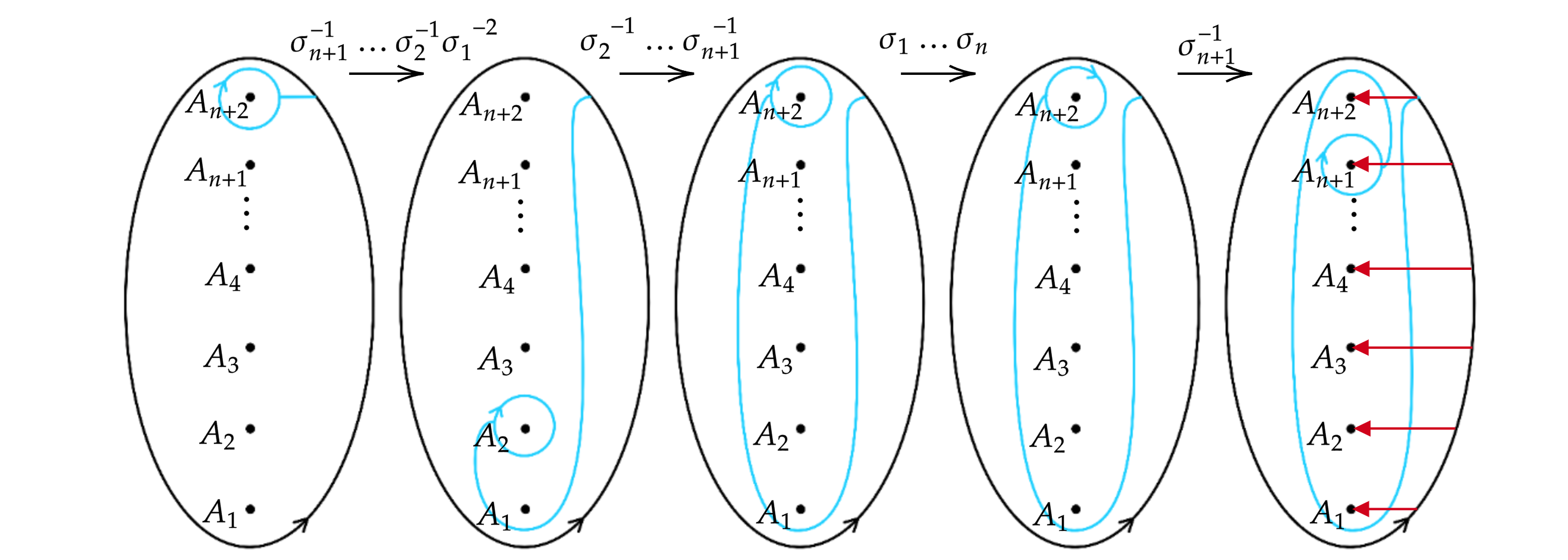}
         \caption{Deduction of the branching word, case 2}
    \label{branching3}
     \end{figure}
Following the Figure \ref{branching3}, associated with the $A_{n+2}$ strand, we have the branching word: 
$$W=A_{n+1}^{-1}A_n^{-1}\dots A_{1}^{-1}A_{n+2}^{-1}A_{n+1}^{-1}A_{n+2}A_1\dots A_nA_{n+1}A_{n+2}$$

Defining $\Omega=A_{n+2}A_1A_2\dots A_{n}$, we obtain $W=A_{n+1}^{-1}\Omega^{-1}A_{n+1}^{-1}\Omega A_{n+1}A_{n+2}$.

 This branching word has a different structure, and admits more simplifications than the previous cases. First, we will prove something similar:

\textbf{Claim:} There exists an $i$ such that:
\begin{equation}\tag{C3}\label{cond3}
\begin{matrix}
    a_{n+1}(i)\neq i\\ a_{n+1}^{-1}a_{n+2}^{-1}a_{1}^{-1}\dots a_{n}^{-1}(i)\neq a_{n+2}^{-1}a_{1}^{-1}\dots a_{n}^{-1}(i) 
\end{matrix}
     \end{equation}
Clearly, we have at least one $i$ that satisfies the first part of this condition, so again, by contradiction, suppose that the second part is not true for any of those $i$. If we write the detailed branching word as:
$$
A_{n+1,ji}^{-1}\Omega^{-1}_{ik}A_{n+1,ks}^{-1}\Omega_{st} A_{n+1,tu}A_{n+2,uj}
$$
The first condition is equivalent to $i\neq j$, and the second one to $k=s$. After reducing this detailed branching word as in the previous case, we obtain:
\begin{align*}
A_{n+1,ji}^{-1}\Omega^{-1}_{is}A_{n+1,ss}^{-1}&\Omega_{st} A_{n+1,tu}A_{n+2,uj}\\
A_{n+1,ji}^{-1}\Omega^{-1}_{is}&\Omega_{st} A_{n+1,tu}A_{n+2,uj}\\
A_{n+1,ji}^{-1} &A_{n+1,iu}A_{n+2,uj}\\
&A_{n+2,jj}\\
&\emptyset
\end{align*}
This means that if $i$ is a non-fixed number in the permutation $a_{n+1}$, then $a_{n+1}(i)$ cannot be in $a_{n+2}$. It follows that $a_{n+1}$ and $a_{n+2}$ have to be disjoint permutations, and they commute.

By using the coloring of Figure \ref{coloration3}, we obtain the Wirtinger relations:
\begin{align*}
    \tag{1} b_{n-1}&=a_{n+1}a_{n+2}a_{n+1}^{-1}\\
    \tag{2} b_i&=a_{i+2}b_{i+1}a_{i+2}^{-1}& 0\leq i\leq n-2\\
    \tag{3} c_{n}&=a_1b_0a_{1}^{-1}\\
    \tag{4} c_0&=c_na_1c_{n}^{-1}\\
    \tag{5} a_{i}&=c_na_{i+1}c_n^{-1} & 1\leq i \leq n\\
    \tag{6}c_i&=a_{i}^{-1}c_{i-1}a_{i} &1\leq i\leq n-1\\
    \tag{7}a_{n+1}&=a_{n+2}c_{n}a_{n+2}^{-1}
\end{align*}

Since $a_{n+1}a_{n+2}=a_{n+2}a_{n+1}$, then:

\begin{align*}
    a_i&=a_{n+1}\hspace{1cm}1\leq i\leq n+1\\
    c_i&=a_{n+1}\hspace{1cm}0\leq i\leq n-1\\
    b_i&=a_{n+1}\hspace{1cm}0\leq i\leq n-1
\end{align*}
Therefore, by $(1)$ and the previous equalities: $a_{n+1}=a_{n+2}$, and again this is a contradiction since it would imply that $a_{n+1},a_{n+2}$ are nontrivial, equal, disjoint permutations.

We now analyze other possible simplifications in the branching word $W$. To simplify notation, and without loss of generality, we assume $j = 1$ and $i =2$ in condition \ref{cond3}. In other words, the detailed branching word is:

\[
A_{n+1,12}^{-1}\Omega^{-1}_{2k}A_{n+2,kl}^{-1}A_{n+1,lm}^{-1}A_{n+2,mx}\Omega_{xy} A_{n+1,yz}A_{n+2,z1}
\]
with $l \neq m$, and where $x = s$, $y = t$, and $z = u$. Denote by $\omega = a_n \dots a_1$ the permutation associated to $\Omega$. Note that since $A_{n+1,12}^{-1}$ and $A_{n+1,lm}^{-1}$ do not simplify on their own, they can only be simplified with each other or with $A_{n+1,yz}$.

In particular, any simplification resulting in the empty word, or a word beginning with $A_{n+2}$, necessarily requires that the terms between $A_{n+1,12}^{-1}$ and $A_{n+1,lm}^{-1}$ (that is, $\Omega^{-1}_{2k}A_{n+2,kl}^{-1}$) must simplify among themselves. This can only occur if all such terms fix the sheet $2$, i.e., the factor must take the form:

\[
\Omega^{-1}_{22}A_{n+2,22}^{-1} = A_{n,22}^{-1}A_{n-1,22}^{-1}\dots A_{1,22}^{-1}A_{n+2,22}^{-1}
\]

In this scenario, we have $a_k(2) = 2$ for $1 \leq k \leq n$ and $k = n+2$. Plugging this into the detailed branching word, we get:

\begin{equation}\label{eq1}
A_{n+1,12}^{-1}\Omega_{22}^{-1}A_{n+2,22}^{-1}A_{n+1,2m}^{-1}A_{n+2,mx}\Omega_{xy}A_{n+1,yz}A_{n+2,z1}
\end{equation}

which reduces to:
\[
A_{n+1,12}^{-1}A_{n+1,2m}^{-1}A_{n+2,mx}\Omega_{xy}A_{n+1,yz}A_{n+2,z1}
\]

We now distinguish two cases:

\begin{itemize}
    \item If $m \neq 1$, we obtain a simplification of type 3:
    \[
    A_{n+1,1m}A_{n+2,mx}\Omega_{xy}A_{n+1,yz}A_{n+2,z1}
    \]
    The only potentially problematic simplification occurs when $m = x = y$:
    \begin{align*}
    A_{n+1,1x}A_{n+2,xx}&\Omega_{xx}A_{n+1,xz}A_{n+2,z1} \\
        A_{n+1,1x}&A_{n+1,xz}A_{n+2,z1}
    \end{align*}
    
    However, since $A_{n+1,1x}$ arises from a type 3 simplification with $m \neq 1$, the terms $A_{n+1,1x}A_{n+1,xz}$ cannot be reduced to the empty word. It follows that the detailed branching word cannot be simplified to the empty word or to one starting with $A_{n+2}$, and thus the monodromy is left-veering.

    \item If $m = 1$, then (\ref{eq1}) becomes:
    \begin{equation}\label{eq2}
    A_{n+1,12}^{-1}\Omega_{22}^{-1}A_{n+2,22}^{-1}A_{n+1,21}^{-1}A_{n+2,1x}\Omega_{xy}A_{n+1,yz}A_{n+2,z1}
    \end{equation}
\begin{figure}
    \centering
    \includegraphics[width=1\linewidth]{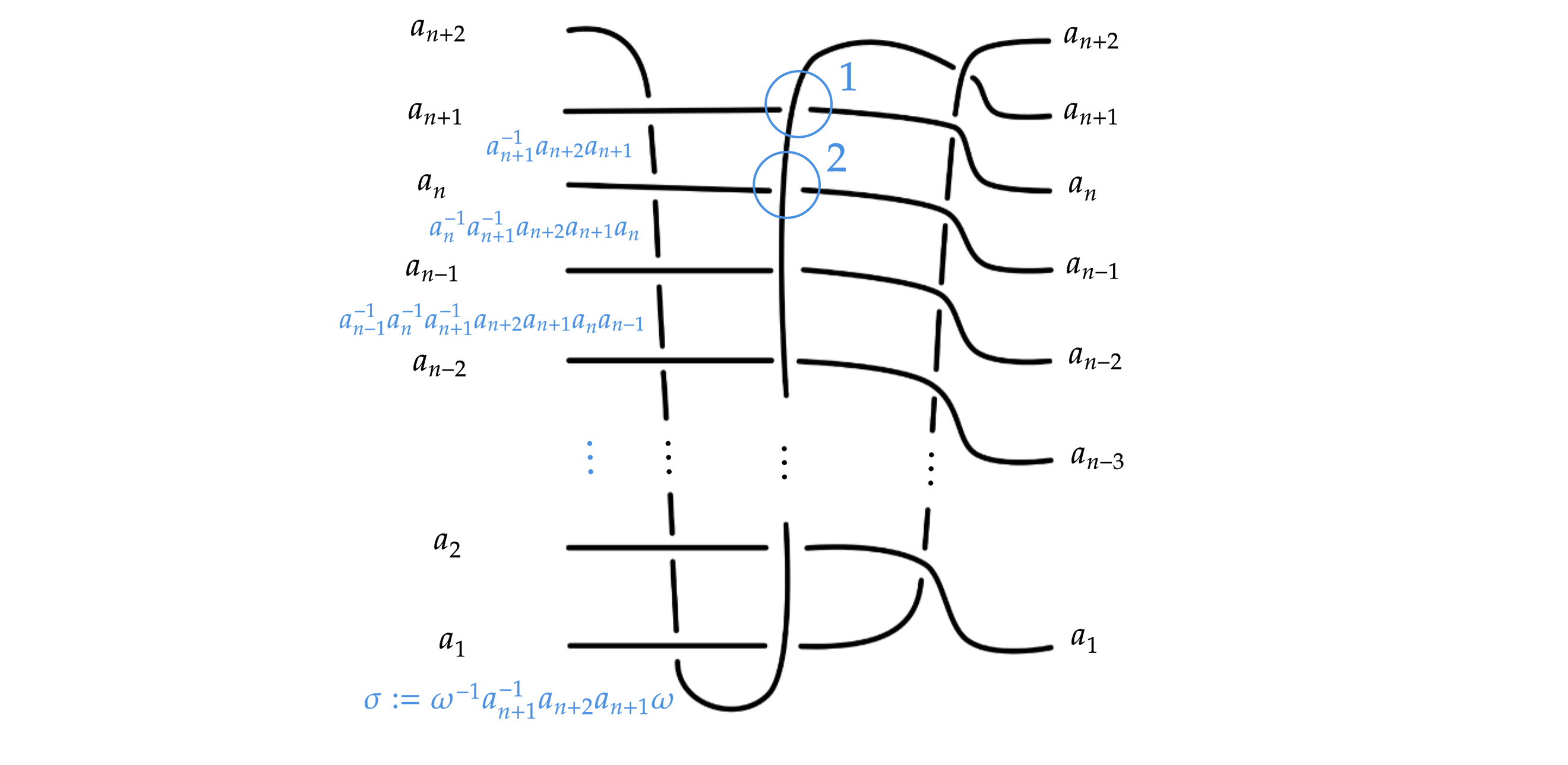}
    \caption{Marked crosses on the braid for $K_m$, $m\geq 2$ odd}
    \label{Two marked crosses}
\end{figure}
    
    In particular, the permutation $a_{n+1}$ contains the transposition $(1\hspace{0.3cm} 2)$. Starting from sheet $z$ in (\ref{eq2}), we can rewrite it as:
    \begin{equation}\label{eq3}
    A_{n+1,zy}^{-1}\Omega_{yx}^{-1}A_{n+2,x1}^{-1}A_{n+1,12}^{-1}A_{n+2,22}\Omega_{22}A_{n+1,21}A_{n+1,1x}
    \end{equation}
    
    This implies $x = z$. Substituting this into (\ref{eq2}) gives:
    \[
A_{n+1,12}^{-1}\Omega_{22}^{-1}A_{n+2,22}^{-1}A_{n+1,21}^{-1}A_{n+2,1x}\Omega_{xy}A_{n+1,yx}A_{n+2,x1}
    \]

    If $x = 1$, then $y = 2$, but this contradicts the fact that $2$ is fixed by $\omega$ since $\omega(x)=y$. Similarly, if $x = 2$, then $y = 1$, which also contradicts the fixed point condition. Therefore, $x \notin \{1,2\}$, and for simplicity, we set $x = 3$. Replacing into (\ref{eq3}):
    \begin{align*}
        A_{n+1,3y}^{-1}\Omega_{y3}^{-1}A_{n+2,31}^{-1} & A_{n+1,12}^{-1}A_{n+2,22}\Omega_{22}A_{n+1,21}A_{n+1,13} \\
        A_{n+1,3y}^{-1}\Omega_{y3}^{-1}&A_{n+2,31}^{-1}  A_{n+1,12}^{-1}A_{n+1,21}A_{n+1,13} \\
         A_{n+1,3y}^{-1} & \Omega^{-1}_{y3}
    \end{align*}
    
    Observe that if $y \neq 3$, then $A_{n+1}^{-1}$ cannot simplify, so the monodromy remains left-veering.

    Suppose $y = 3$. Then the permutations partially take the form:
    \begin{align*}
        a_k &= (2)\dots\hspace{2cm}\text{for } 1 \leq k \leq n, \\
        \omega &= (2)(3)(1\hspace{0.3cm} p\hspace{0.3cm} \dots)\dots, \\
        a_{n+1} &= (1\hspace{0.3cm} 2)(3)\dots, \\
        a_{n+2} &= (1\hspace{0.3cm} 3)(2)\dots,
    \end{align*}
    where the omitted terms involve none of $1,2,$ or $3$, and $p\notin\{1,2,3\}$ is some natural number.

Now consider again the Wirtinger relations for the braid in Figure \ref{Two marked crosses}. From the marked crossings 1 and 2, we obtain:
\[
a_n = \sigma^{-1} a_{n+1} \sigma, \qquad
a_{n-1} = \sigma^{-1} a_n \sigma = \sigma^{-2} a_{n+1} \sigma^2
\]
where $\sigma = \omega^{-1}a_{n+1}^{-1}a_{n+2}a_{n+1}\omega$.

Using subscripts to trace the permutations (e.g., $a_k(1) = 3$ is written as $a_{k,31}$), and since $a_{n-1}(2) = 2$, we read the equality $a_{n-1} = \sigma^{-2}a_{n+1}\sigma^2$ from right to left as:

\[
a_{n-1,22} = \omega^{-1}a_{n+1}^{-1}a_{n+2}^{-2}a_{n+1}\omega_{p1} a_{n+1,12} \omega_{22}^{-1}a_{n+1,21}^{-1}a_{n+2,11}^2 a_{n+1,12}\omega_{22},
\]
and from left to right as:
\[
a_{n-1,22} = \omega_{22}^{-1}a_{n+1,21}^{-1}a_{n+2,11}^{-2}a_{n+1,12}\omega_{22} a_{n+1}\omega^{-1}a_{n+1}^{-1}a_{n+2}^2 a_{n+1} \omega.
\]

Comparing subscripts from the fifth term $\omega$ yields $p = 2$ and $1 = 2$, each of which results in a contradiction. Hence, this case is not possible.
\end{itemize}

It follows that in all cases, we can always find a non-empty reduced detailed branching word that does not begin with $A_{n+2}$, and thus the monodromy is left-veering.

\begin{corollary} For $m\geq2$ and $m=2k-1\leq-3$, every transverse representative of $K_m$ is not transverse-universal.
\end{corollary}

\section{Case \texorpdfstring{$m=-2n\leq-4$}{m=-2n<=-4}}

Theorem \ref{classification} implies that, except from the case $m=-4$, we do not have transverse simplicity for this subfamily. We could think of doing the same analysis as before for each possible transverse representative with maximal self-linking number, proving in that way that they also have only overtwisted branched covering spaces. However, this is far from the truth. We will see that each of these transverse knots has infinitely many associated tight branched coverings.

The first problem is that we need to be able to distinguish different transverse representatives of maximal self-linking number that belong to the same topological class. Fortunately, Etnyre, Ng, and V\'ertesi gave an answer for Legendrian representatives and therefore, for transversal as well; see \cite{EtLenVer} for further details. 

\begin{figure}
    \centering
    \includegraphics[width=14cm]{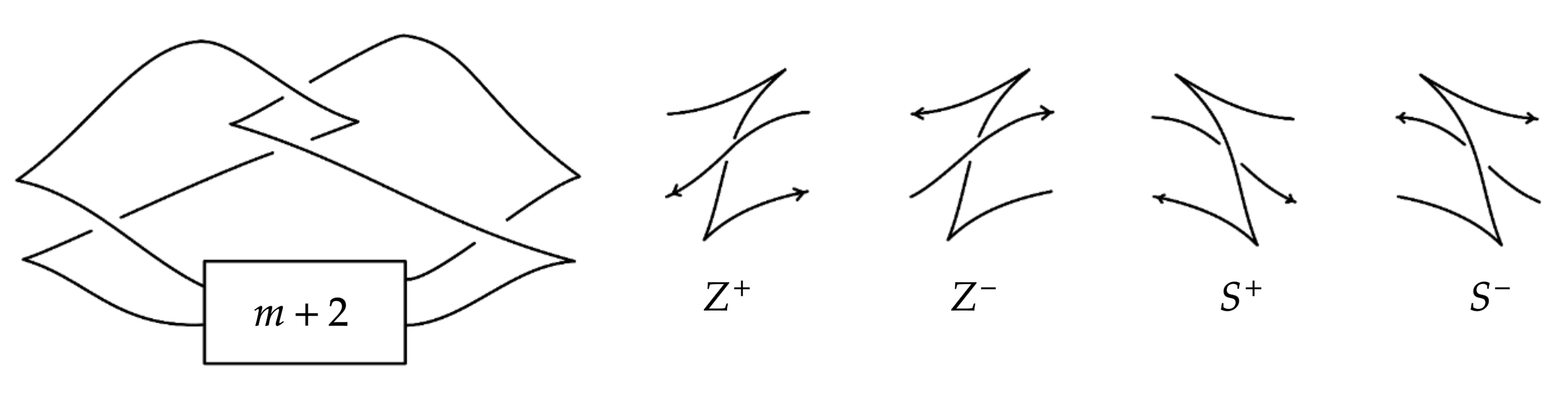}
    \caption{Codification for legendrian representatives}
    \label{CodLeg}
\end{figure}

\begin{theorem}[\cite{EtLenVer}]\label{ClasLeg}
For $m=-2n\leq -2$, any Legendrian representative of $K_m$ with maximal $tb$ is Legendrian isotopic to some legendrian knot whose front projection is of the form depicted in Figure \ref{CodLeg}, where each of the $m+2$ crossings is a $Z^{\pm}$ or $S^{\pm}$.
    
\end{theorem}
\begin{figure}
    \centering
    \includegraphics[width=14cm]{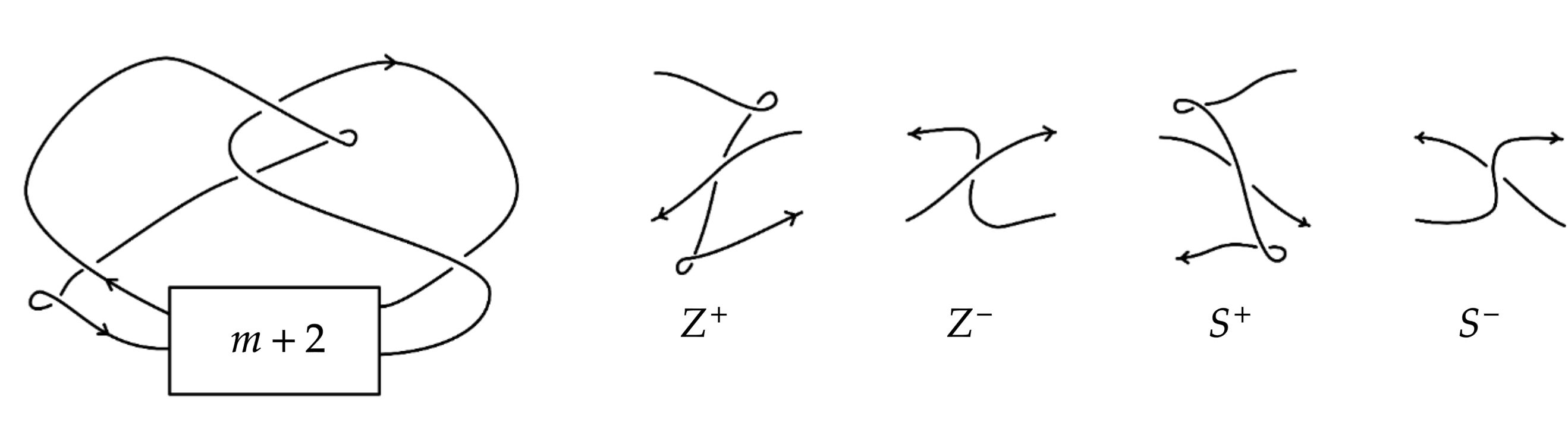}
    \caption{Codification for transverse representative}
    \label{CodTran}
\end{figure}

Denote by $z^\pm$, $s^\pm$ the number of crossings of type $Z^\pm$ or $S^\pm$ respectively. It is easy to note that $z^++s^+=z^-+s^-=n-1$, since each positive crossing ($Z^+,S^+$) can only be next to a negative one ($Z^-,S^-$), as the figure shows. Therefore, if we fix $m$ and, for example, $z^+,z^-$, we can recover all the Legendrian representatives with maximal $tb$ for those values. A Legendrian representative of $K_m$, like the previous ones, will be denoted by $K_m(z^+,z^-)$.

\begin{proposition}[\cite{EtLenVer}]\label{ClasLeg2}
    Let $m=-2n\leq -4$, then $K_m(z^+,z^-), K_m(\tilde z^+,\tilde z^-)$ are Legendrian isotopic after some positive number of negative stabilizations if and only if $\tilde z^+=z^+$ or $\tilde z^+=n-1-z^+$, and in these cases the knots are isotopic after one negative stabilization. 
\end{proposition}

\begin{figure}
    \centering
    \includegraphics[width=11cm]{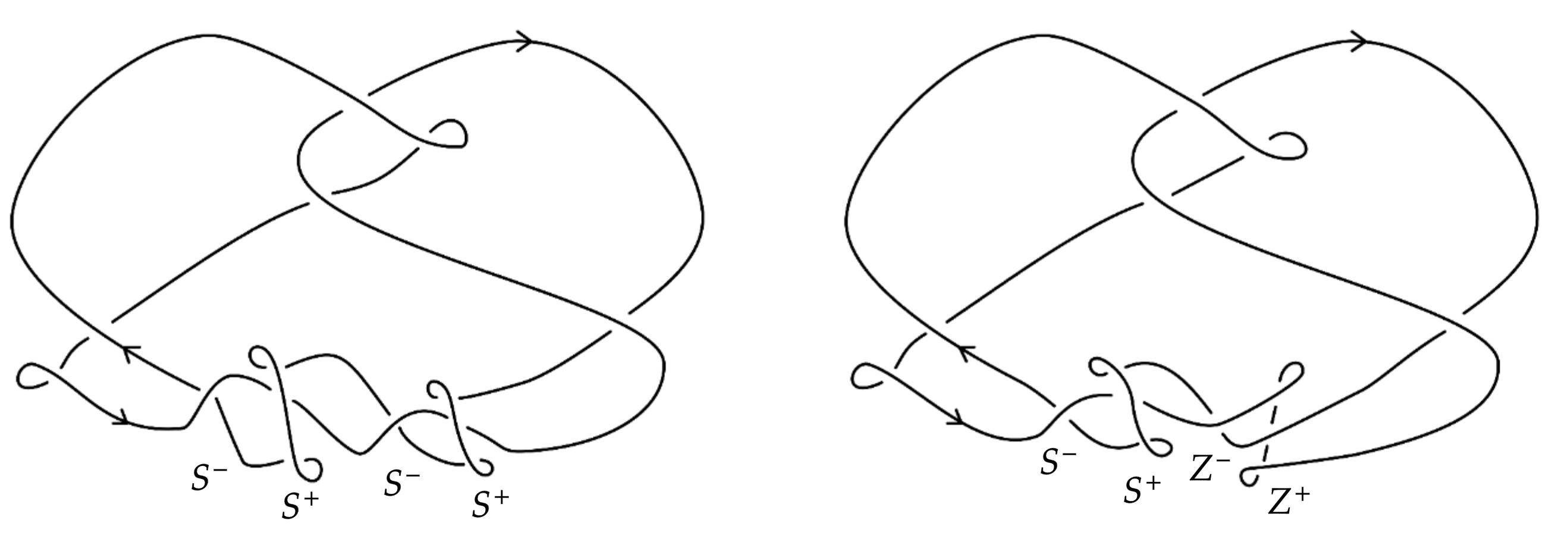}
    \caption{Example for $m=-6$}
    \label{CodExam}
\end{figure}

We also know that for $m=-2n\leq-4$ even, $K_m$ has $\lceil \frac{n^2}{2} \rceil$ different Legendrian representatives with $(tb,rot)=(1,0)$, and this is the maximal value for the $tb$ of this subfamily \cite{EtLenVer}. Even more, a direct application of  Proposition \ref{ClasLeg2} implies that these Legendrian knots fall into $\lceil \frac{n}{2} \rceil$ different Legendrian isotopy classes after any given positive number of negative stabilizations.

In this way, we can use the codification for Legendrian representatives and consider their transverse push-offs with orientation as in Figure \ref{CodTran}. For example, Figure \ref{CodExam} shows the case $m=-6$ with two different values of $z^+,z^-$, and since $\lceil \frac{3}{2} \rceil=2$, we can conclude that these two transverse knots are all the transverse representatives of $K_{-6}$ with maximal self-linking number $1$. 

In general, for each $1\leq l\leq\lceil \frac{n}{2} \rceil $, we have a transverse representative of $K_m=K_{-2n}$, and all of them are listed in the following table with the associated values of $z^+,s^+,z^-,s^-$:

{\centering
\begin{tabular}{|c|c|c|c|c|c|c|}
\hline
     &1&2&3&4&...& $\lceil \frac{n}{2} \rceil$ \\\hline
     $z^+$&0&1&2&3&...&$\lceil \frac{n}{2} \rceil-1$\\\hline
     $s^+$&$n-1$&$n-2$&$n-3$&$n-4$&...&$\lceil \frac{n}{2} \rceil$\\\hline
     $z^-$&0&1&2&3&...&$\lceil \frac{n}{2} \rceil-1$\\\hline
     $s^-$&$n-1$&$n-2$&$n-3$&$n-4$&...&$\lceil \frac{n}{2} \rceil$\\\hline
\end{tabular}
\par}

\begin{figure}
    \centering
    \includegraphics[width=14cm]{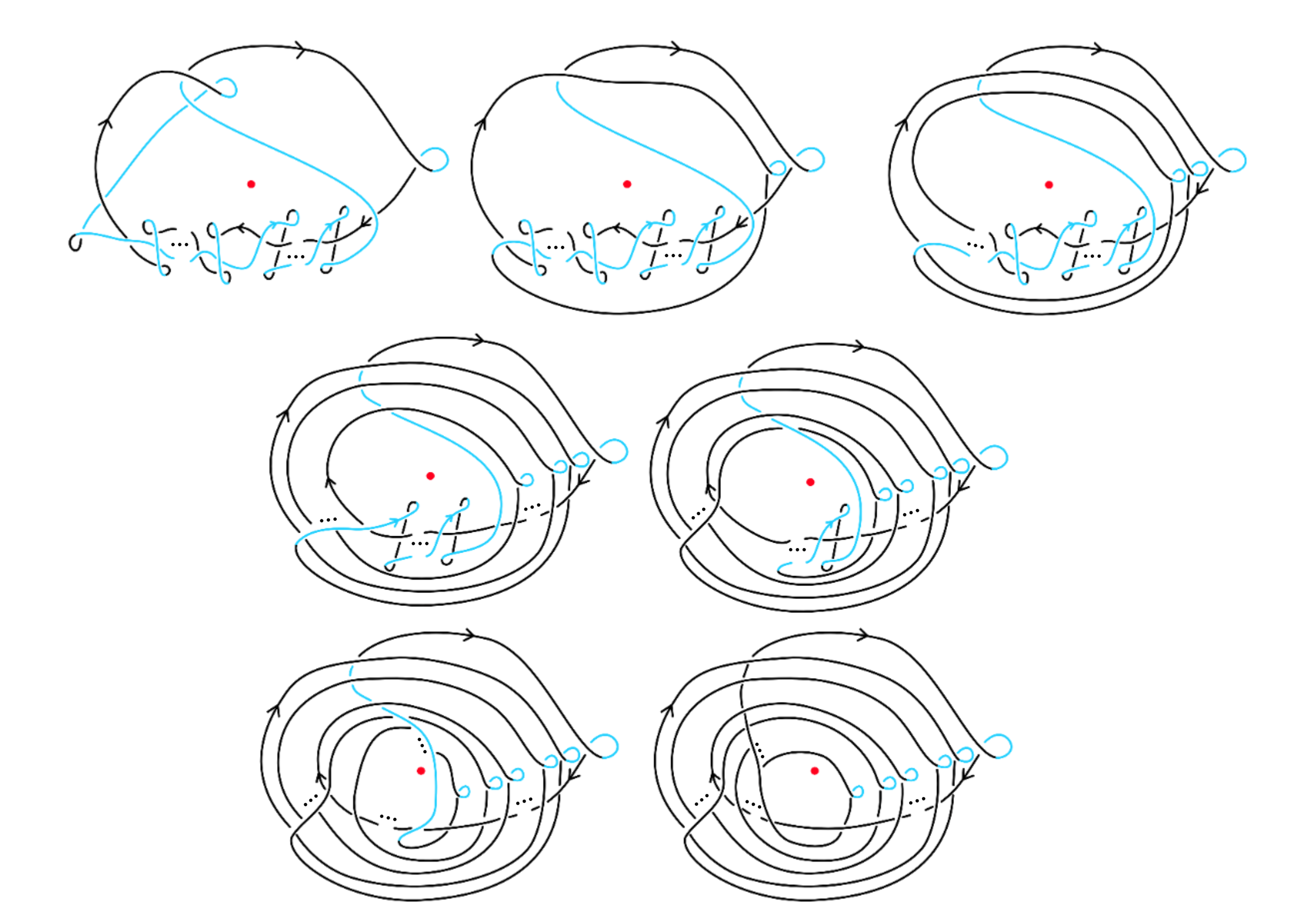}
    \caption{Braiding process}
    \label{BraidingEven}
\end{figure}
For each column, we choose the order
$$S^-S^+S^-\dots S^-S^+Z^-Z^+\dots Z^-Z^+$$

and hence, we have a preferred way to represent the desired members of this subfamily.

\begin{figure}
\centering
\begin{minipage}{.5\textwidth}
  \centering
    \includegraphics[width=9.5cm]{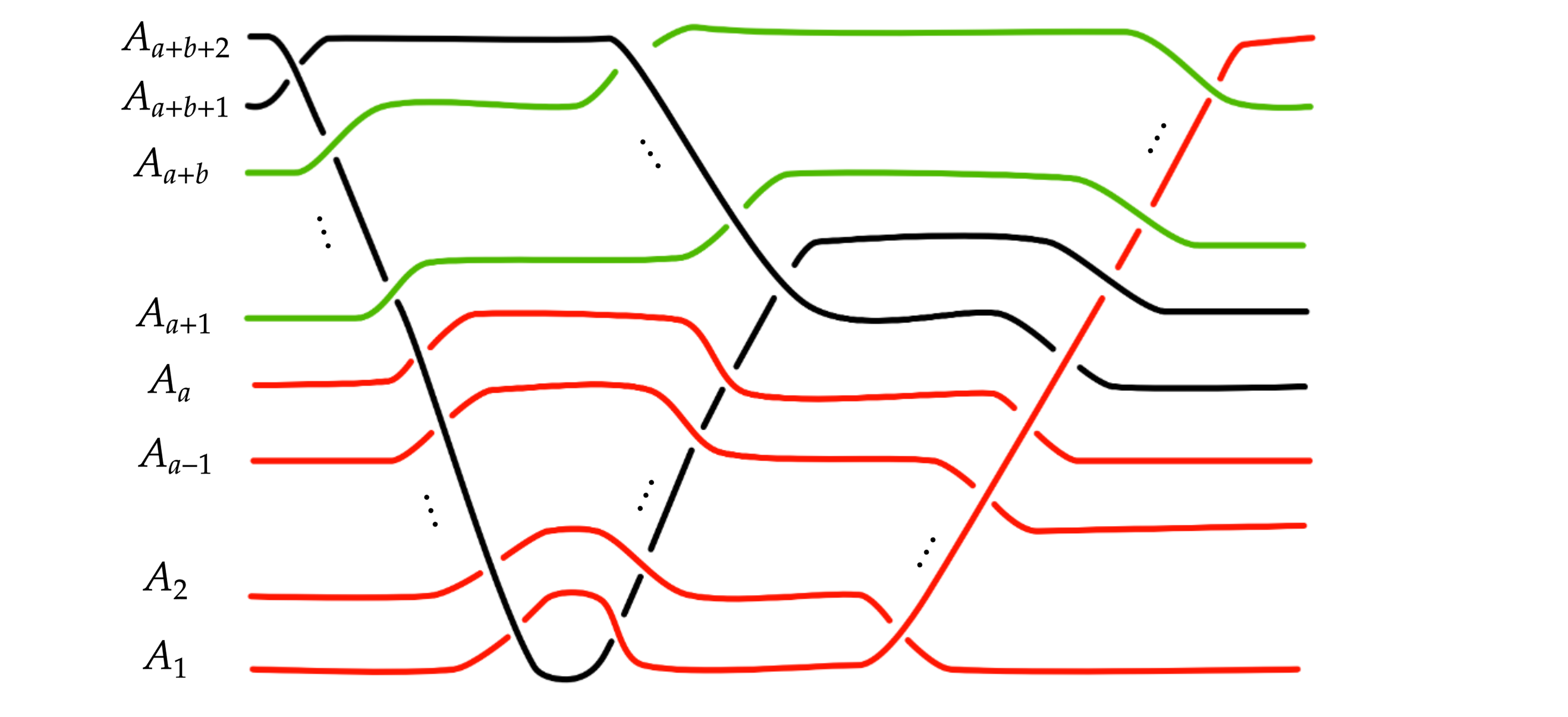}
    \caption{Case $a>0$}
    \label{FirstCase}
\end{minipage}%
\begin{minipage}{.5\textwidth}
  \centering
    \includegraphics[width=9.5cm]{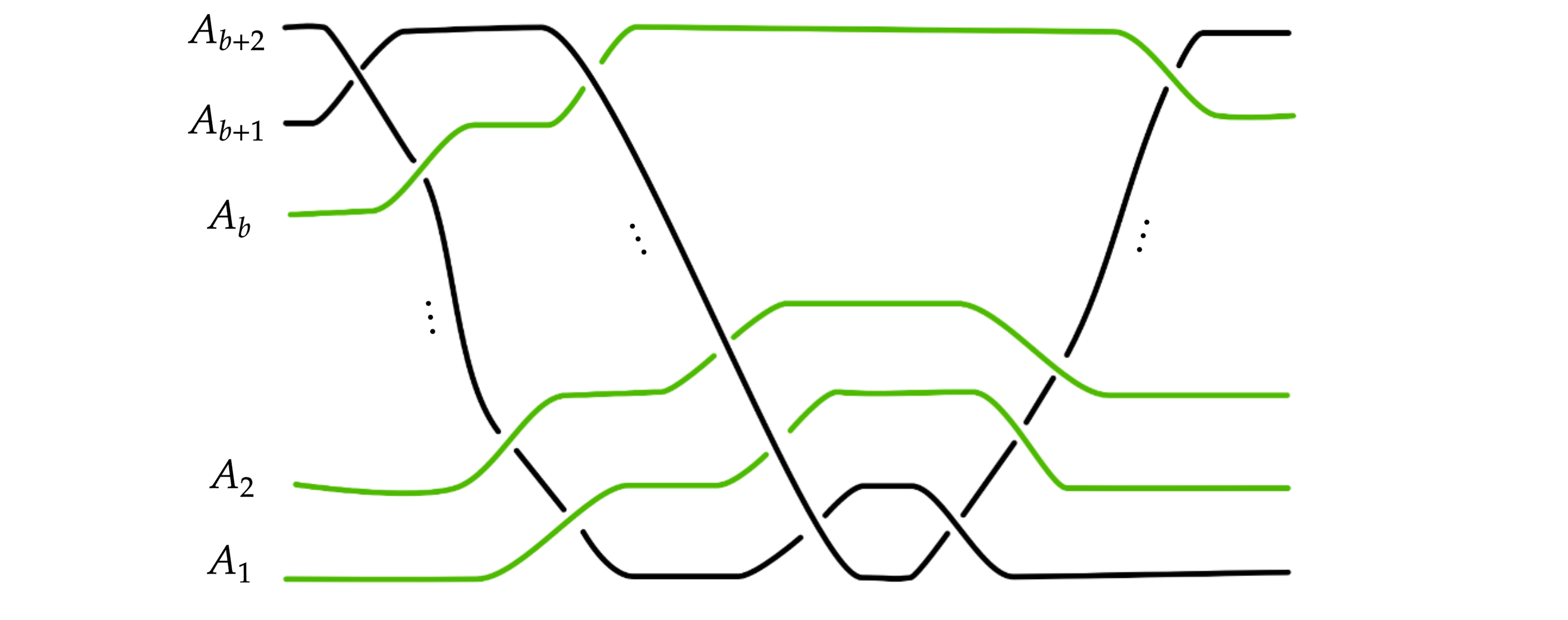}
    \caption{Case $a=0$}
    \label{SecondCase}
\end{minipage}
\end{figure}

Let $a=z^+$, and $b=s^+$. The general braiding process is shown in Figure \ref{BraidingEven}, but to analyze the braid, we need to distinguish two possibilities:

\begin{itemize}
    \item [(1)] $a>0:$ In this case, we obtain the braid in Figure \ref{FirstCase}. The green and red strands are induced by the $S^\pm, Z^\pm$ crossings, respectively, and we have the braid words:
    
For $a\geq2:$
$$
\sigma_{a+b+1}\sigma_{a+b}^{-1}\dots \sigma_{a+1}^{-1}\sigma_{a}\dots \sigma_1\sigma_1\dots\sigma_a\sigma_{a+b+1}\dots\sigma_{a+1}\sigma_{1}^{-1}\dots\sigma_{a}^{-1}\sigma_{a+1}\dots\sigma_{a+b+1}
$$
and this is equivalent to
$$
\sigma_{a+b+1}[\sigma_1\dots\sigma_a\sigma_{a+b+1}\dots\sigma_{a+1}]^{\sigma_1^{-1}\dots\sigma_{a}^{-1}\sigma_{a+1}\dots\sigma_{a+b}}\sigma_{a+b+1}
$$
For $a=1:$
$$
\sigma_{b+2}\sigma_{b+1}^{-1}\dots \sigma_{2}^{-1}\sigma_{1}\sigma_1\sigma_{b+2}\dots\sigma_{2}\sigma_{1}^{-1}\sigma_{2}\dots\sigma_{b+2}
$$

$$
\sigma_{b+2}[\sigma_{1}\sigma_{b+2}\dots\sigma_{2}]^{\sigma_{1}^{-1}\sigma_2\dots\sigma_{b+1}}\sigma_{b+2}
$$

\item [(2)] $a=0:$ In this case, the braid takes the form of Figure \ref{SecondCase}, and the braid word is:
$$
\sigma_{b+1}\sigma_{b}^{-1}\dots\sigma_{1}^{-1}\sigma_{b+1}\dots\sigma_{1}\sigma_{1}\dots\sigma_{b+1}
$$
$$
\sigma_{b+1}[\sigma_{b+1}\dots\sigma_{1}]^{\sigma_{1}\dots\sigma_{b}}\sigma_{b+1}
$$
\end{itemize}

In each case, we have that the braid is quasipositive; hence, every cyclic contact branched covering, with branch locus one of these transverse knots is Stein fillable, and thus tight \cite[Theorem 1.3]{HaKaPla}. The specific case of double coverings was observed by Plamenevskaya in \cite{Pla}.

 \textbf{Remark:} The only twist knots not covered so far correspond to $m=-2,-1,0,1$. The knots $K_{-1}$ and $K_0$ are both trivial, while $K_{-2}$ and $K_{1}$ are the positive and negative trefoil knots, respectively. None of these is topologically universal, but one can verify that they admit tight contact branched coverings. For example, a cyclic branched covering of $K_{-2}$ and the dihedral triple covering of $K_1$

\bibliographystyle{plain}
\bibliography{Bib}
\end{document}